# Lectures on Mean Curvature Flow
# and Related Equations

*Draft Version, recompiled August 1998*

## Tom Ilmanen

Former title: "Some Techniques in Geometric Heat Flow." These notes were prepared in $\mathcal{AMS}$-TeX for the Conference on Partial Differential Equations & Applications to Geometry, 21 August - 1 September, 1995, ICTP, Trieste. The author was supported in part by an NSF Postdoctoral Fellowship.

# Lecture 1

Our aim in these lectures is to study singularity formation, nonuniqueness, and topological change in motion by mean curvature. The outline is:

Lecture 1 - Introduction

Lecture 2 - Flows in $\mathbf{R}^3$

Lecture 3 - Flows in $\mathbf{R}^3$, continued

Lecture 4 - Nonuniqueness in Geometric Heat Flows

Let $\{M_t\}_{t \in \mathbf{R}}$ be an evolving family of hypersurfaces in $\mathbf{R}^n$. We say that $M_t$ is moving by mean curvature if it satisfies the nonlinear parabolic equation

(MCF) $$\frac{\partial}{\partial t} x = \vec{H}(x), \qquad x \in M_t, \quad t \in \mathbf{R}.$$

Here $\vec{H}(x)$ is the mean curvature vector of $M_t$ at $x$, defined to be $\sum_1^{n-1} \lambda_i \nu$, where $\lambda_i$ are the principal curvatures, and $\nu$ is the unit normal. The expression $\partial x/\partial t$ stands for the normal velocity of the surface.

Motion by mean curvature is the gradient flow of the area functional (with respect to the $L^2$ norm on the surface). We have the following formula for the decrease of the area of $M_t$,

$$\frac{d}{dt} \mathcal{H}^{n-1}(M_t) = -\int_{M_t} H^2 \, d\mathcal{H}^{n-1}.$$

Mean curvature flow arises as a simplified model for various physical processes in which surface tension plays a role.

Since (MCF) is a parabolic equation, it will have a short-term smoothing effect on the surface. However, the surface can become singular later.

The simplest example is a sphere. The sphere shrinks without changing shape since the curvature is the same all around. We find that (MCF) becomes

$$\frac{dR}{dt} = -\frac{n-1}{R}, \qquad R(0) = R_0,$$

where $R = R(t)$ is the radius of the sphere, which yields

$$R(t) = \sqrt{R_0^2 - 2(n-1)t},$$

so the sphere disappears at time $R_0^2/2(n-1)$.

Let me mention another early result. Huisken [Hu1] proved that if $M_0$ is a bounded convex surface, then $M_t$ becomes more and more nearly spherical as it shrinks, and at the instant it vanishes, it is asymptotic to the shrinking sphere given above. (See also Hamilton [H8].)

If we could control the curvature of $M_t$, we could keep it from becoming singular. But by differentiating (MCF), we have the following equation for the evolution of the norm of the extrinsic curvature (see [Hu1]):

$$\frac{\partial}{\partial t}|A|^2 = \Delta_{M_t}|A|^2 - 2|D_{M_t}A|^2 + 2|A|^4.$$

Here $A = (A_{ij})_{i,j=1}^{n-1}$ is the second fundamental form of $M_t$ in $\mathbf{R}^n$, $|A|$ is its length, and $\Delta_{M_t}$, $D_{M_t}$ are intrinsic derivatives on the surface $M_t$. Note that $\partial|A|^2/\partial t$ means the derivative with respect to a "meterial point" moving perpendicularly to the surface. This equation resembles

$$u_t = \Delta u + u^2$$

which has solutions that blow up to infinity. This suggests that it will be difficult to control the singularities. This is borne out by the example of the sphere, where the motion accelerates as the sphere gets smaller.

Generally speaking, the singularities of motion by mean curvature can be very complicated.

The mean curvature flow shares many characteristics with other so-called *geometric heat flows*, such as the harmonic map flow, the equation $u_t = \Delta u + u^p$, the Yang-Mills heat flow, and the Ricci flow. In particular, many of the ideas presented in these lectures will apply to these other geometric heat flows. I will have more to say about this in the last lecture.

**Main Questions.** Here are some major issues in mean curvature flow; they are at various stages of development by a number of mathematicians.

*1. Nature of singularity formation, especially in $\mathbf{R}^3$*

The basic mechanism is self-similar shrinking, as in the case of the convex bodies becoming spherical, but there are other, more subtle mechanisms, see Angenent-Velazquez [AV1, AV2], Velazquez [V].

## 2. *Weak solutions past the onset of singularities*

Two definitions of weak solution have gained attention in recent years:

- the *level-set flow*, based on the idea of viscosity solutions (tangency with smooth test functions), see [ES1-ES4], [CGG],

- *Brakke's flow*, using the varifolds of geometric measure theory (integration by parts with the surface measure), see [B].

## 3. *Partial regularity*

The best result so far is : almost every evolution is smooth almost everywhere [I3]. It is based on Brakke's local regularity theorem [B], together with the level-set flow. Recent work of White [W2-W4,W6] points towards a deeper theory of partial regularity.

## 4. *Selection principle, in case of nonunique evolution*

As we shall see, the mean curvature flow supports nonunique evolution (past the first singularity). Which evolution is better? Various minimization principles are suggestive [Am], [I3, I5]. Perhaps the selected evolutions would decrease area more rapidly and be more regular.

## 5. *Role of random processes*

Recently there has been progress in making rigorous the stochastic origins of mean curvature flow, that is, various interacting particle systems related to the Ising model, see Katsoulakis-Souganidis [ ], Griffeath [ ], and others. Could this yield clues about which flow to select?

# Lecture 2: Flows in $\mathbf{R}^3$

The outline of this lecture is

      A - Self-Shrinking Surfaces

      B - Examples of Self-Shrinkers in $\mathbf{R}^3$

      C - Nonuniqueness in $\mathbf{R}^3$

      D - Level-Set Flow

      E - Topological Changes Past a Singularity

      F - Evolution of Cones in $\mathbf{R}^3$

We state the following principle, discovered by Huisken [Hu2]:

*Singularity formation is modelled by self-shrinking surfaces.*

This is already exemplified by the above-mentioned evolution of convex hypersurfaces. We will accept this principle without proof for the moment; but let us explain what we mean by shrinking self-similarly.

## A. Self-Shrinking Surfaces

By counting units on either side of the equation, the reader will agree that (MCF) is invariant under *parabolic rescaling*

$$x \to \lambda x, \qquad t \to \lambda^2 t.$$

We look for solutions that are invariant under the same scaling, that is, we impose the ansatz

$$N_t = \sqrt{-t} \cdot N, \qquad t < 0.$$

where $\sqrt{-t} \cdot$ represents a homothety.

Another way to motivate this ansatz is to insert the general expression $a(t) \cdot N_t$ into (MCF) to obtain separate equations for $a$ and $N$. We solve for $a$ to obtain $\sqrt{-t}$ (there are other interesting solutions, which we ignore for now). We also obtain the following elliptic parametric equation for $N$:

(SS) $$H + \frac{x \cdot \nu}{2} = 0, \qquad x \in N.$$

Here $H = \vec{H} \cdot \nu$. This says:

If the surface $N$ begins to move by its mean curvature vector, this produces the same normal velocity as moving by the vector $-x/2$.

That is to say, the surface is shrinking by homothety. We call $N$, and $N_t$, a *self-shrinking surface*.

It turns out that (SS) is the Euler-Lagrange equation of a functional, discovered by Huisken [Hu2]. To find this, let us try a functional of the general form

$$\int_N f(x)\,d\mathcal{H}^{n-1}(x),$$

where $f$ is to be determined. Let $X$ be a $C^1$ vectorfield of compact support. Vary $N$ in the direction of $X$ to yield a family $N_s$ defined by $N^s = \Phi^s(N)$, where $\Phi^s(x) = x + sX(x)$, $s \in (-\varepsilon, \varepsilon)$ is a $C^2$ family of diffeomorphisms. Let $Df$ denote ordinary derivation in the ambient space, let $Df^T$ be the orthogonal projection of $Df$ onto $T_x N$, and let $Df^\perp = Df - Df^T$. Write $\operatorname{div}_N X$ for the divergence of $X$ along $N$, that is, $\sum_{i=1}^{n-1} D_{e_i} X \cdot e_i$ where $e_1 \dots, e_{n-1}$ is an orthonormal basis for $T_x N$. We calculate the first variation in the direction $X$:

$$
\begin{aligned}
\frac{d}{ds}\Big|_{s=0}\int_{N^s} f(x) &= \frac{d}{ds}\Big|_{s=0}\int_N f(\Phi^s(x))|J_{T_x N}\Phi^s(x)| \\
&= \int_N Df \cdot X + f\operatorname{div}_N X \\
&= \int_N Df \cdot X + f\operatorname{div}_N(fX) - Df^T \cdot X \\
&= \int_N -fX \cdot \vec{H} + Df^\perp \cdot X
\end{aligned}
$$

(1)

This is called the *weighted first variation formula*. In the last line we used the usual first variation formula

$$(2) \qquad \int_N \operatorname{div}_N X = \int_N -X \cdot \vec{H} \qquad\qquad \left(= \frac{d}{ds}\Big|_{s=0}\mathcal{H}^{n-1}(N^s)\right)$$

which the reader will recognize as an integration by parts formula. From (1) we obtain the Euler Lagrange equation

$$fH - Df \cdot \nu = 0,$$

which becomes (SS) if we choose $f = e^{-|x|^2/4}$. Accordingly, we define

$$\mathbf{J}[N] = \int_N e^{-|x|^2/4} \, d\mathcal{H}^{n-1}(x).$$

## B. Example of Self-Shrinkers in $\mathbf{R}^3$.

We now present some further solutions of (SS).

A well known self-shrinker is the cylinder $S^1 \times \mathbf{R}$ in $\mathbf{R}^3$.

Figure. Cylinder

There is also a rotationally symmetric torus discovered by Angenent [A5] by ODE methods.

Finding further solutions of (SS) has been difficult since they are always unstable as critical points of $\mathbf{J}$. All the other examples are recent and are based on numerical computations, without proof. The first algorithms and most of the examples are due to Chopp [C].

The idea is to start with a surface in $\mathbf{R}^3$ such as the following, evolve it by mean curvature, and rescale the singularities to find interesting examples of self-shrinkers – we hope.

Figure. A Surface

Typically, the singularities that form will consist of little necks pinching off; these are aymptotic to shrinking cylinders. The reason for this is that one shrinking handle gets a little smaller than the others, and then, because of the high curvature around the handle, the shrinking accelerates and outpaces the other handles.

The more handles, the more ways the necks can get out of balance, and the more unstable the self-shrinking surface will be, if it is ever found.

The easiest way to control the instabilities is to impose a group symmetry, which forces certain handles to shrink off simultaneously. Then we expect to see a new type of singularity.

For example, consider the eight-element group $G = D_2 \ltimes Z_2$, where $D_2$ is the dihedral group of order 4, and $\ltimes$ indicates a semidirect product. $G$ acts on $\mathbf{R}^3$ via the rigid motions

$$(x, y, z) \rightarrow (-x, y, z), \qquad (x, y, z) \rightarrow (x, -y, z); \qquad (x, y, z) \rightarrow (-y, x, -z).$$

Consider the initial surface depicted in the figure, which is invariant under $G$.

Figure. Initial Surface

This topology was suggested by Grayson and Ilmanen. The group forces the two interlocked handles to pinch off simultaneously, if they pinch off at all. Evolving this on the computer, Chopp [C] finds that they do pinch off, and rescaling the surface so that the handle size remains constant, it converges to the surface $N$ shown in the figure.

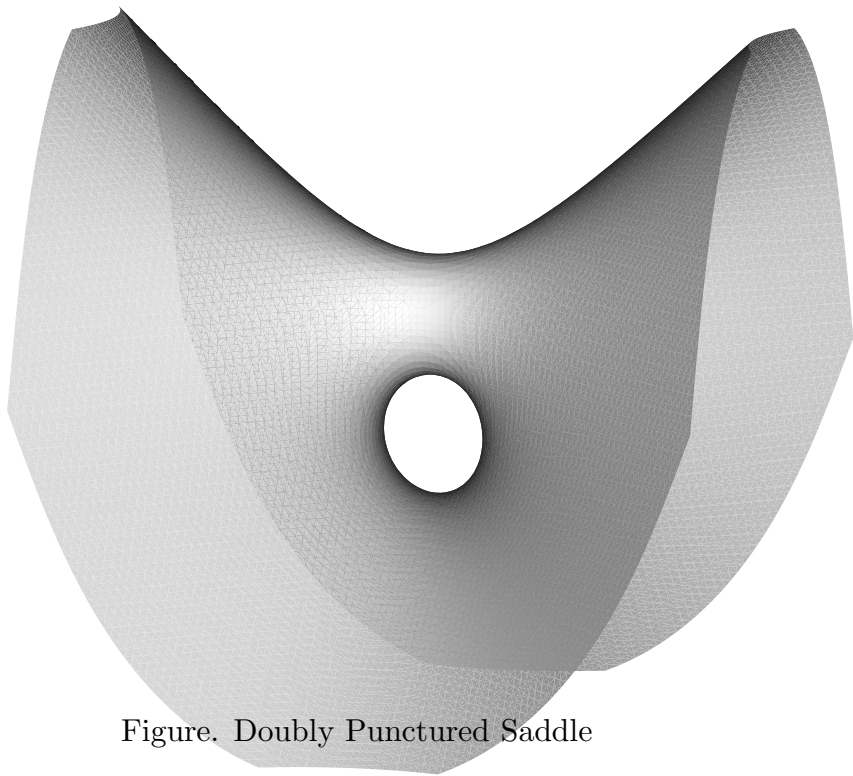

Figure. Doubly Punctured Saddle

Note that $N$ is asymptotic to a certain cone (i.e. homothety invariant set) at spatial infinity. Equivalently, $N_t \equiv \sqrt{-t} \cdot N$ converges to a cone $N_0$ as $t \to 0^-$. We can think of the flow as pulling the surface toward the origin so fast that it streaks out and becomes a cone at the singular time.[1]

We can look for analogous surfaces based on the group $D_m \rtimes Z_2$ for any $m \geq 2$. In the figure is a self-shrinker with $m = 3$. I found it using Brakke's Evolver at the GANG lab at the University of Massachusetts, with technical advice from K. Brakke and J. Sullivan.

I also tried $m = 9$. The resulting surface is depicted in the figure. It looks like the union of a sphere and a plane (each a self-shrinker in its own right) desingularized along the curve of intersection by a high-frequency "zipper" of perforations that looks like the second Scherk surface (the so-called Scherk tower).

These surfaces are the only ones known (besides the sphere and cylinder) where the group symmetry by itself is enough to control the instabilities.

Figure. Scherk Surface

Now we pass to examples which require a mountain-pass argument in addition to the group symmetry. We give the example from Angenent-Chopp-Ilmanen [ACI], using the algorithm from Chopp's paper [C]. Consider the initial surface depicted below, which is a cylinder with holes of a smaller radius bored through it. It is symmetric under the action of the group $G = D_4 \times Z_2$, where the eight-element dihedral group $D_4$ acts on the $y$-$z$ plane (without affecting $x$) and $Z_2$ acts by $(x, y, z) \to (-x, y, z)$.

Because of the group symmetry, the holes pinch off symmetrically, but this time there are two topologically distinct extremes.

---

[1] It can be proven in full generality that as $t \to 0^-$, any self-shrinker converges to a (unique) cone locally in the Hausdorff metric on closed sets. The proof uses the monotonicity formula of Lecture 3, section G; the proof is beyond the scope of the lectures.

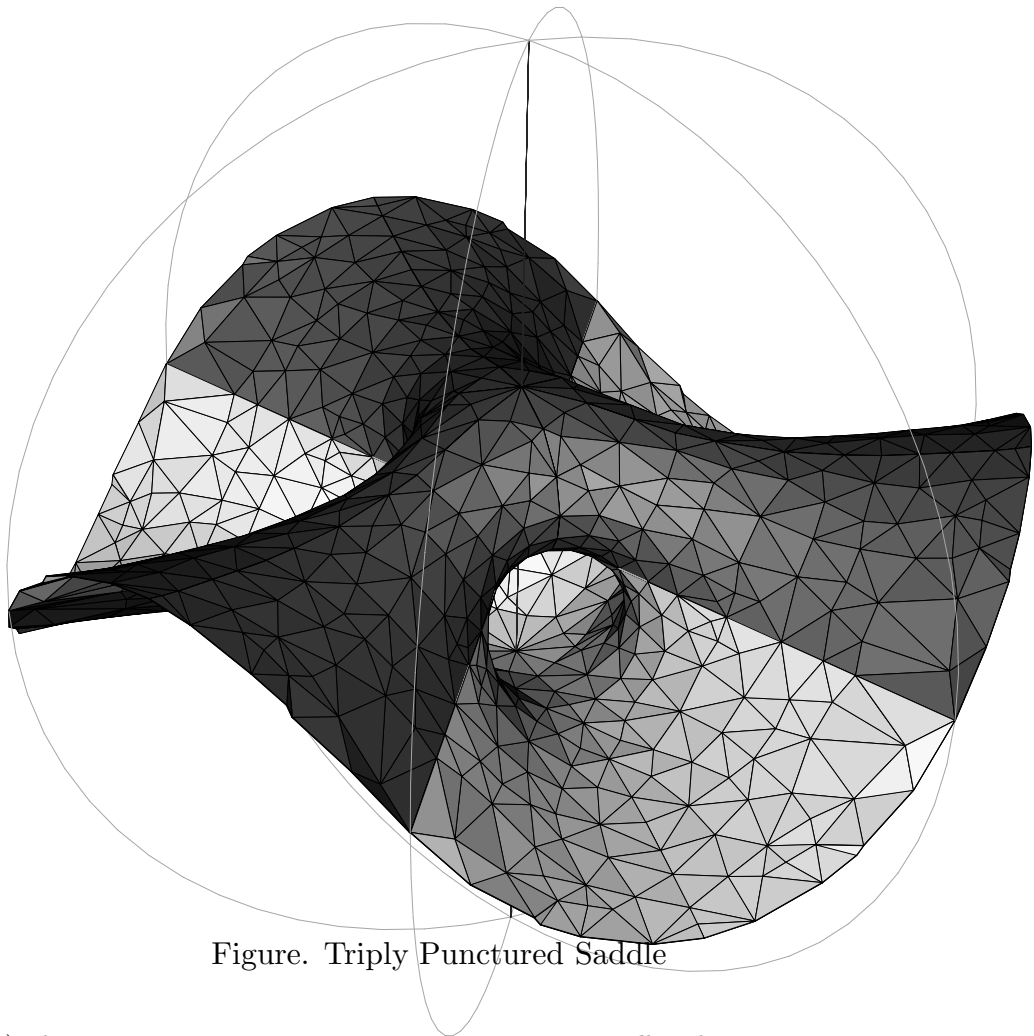

Figure. Triply Punctured Saddle

(a) If the bored holes are very thin, they pinch off in four necks radiating from the origin like spokes of a wheel, leaving a pointy sphere surrounded by a dented cylinder. (After the sphere disappears, the cylinder develops a single neckpinch at the origin.)

(b) If the bored holes are very broad, then the center of the cylinder breaks in four longitudinal necks, with axes parallel to the axis of the cylinder, leaving two pointy stumps, which move away from each other.

Between these two extremes, some third behavior must occur. To find it, we foliated $\mathbf{R}^3$ with the level sets of a $G$-invariant function $u$ whose zero-set is the above surface, and evolved these level-sets simultaneously.

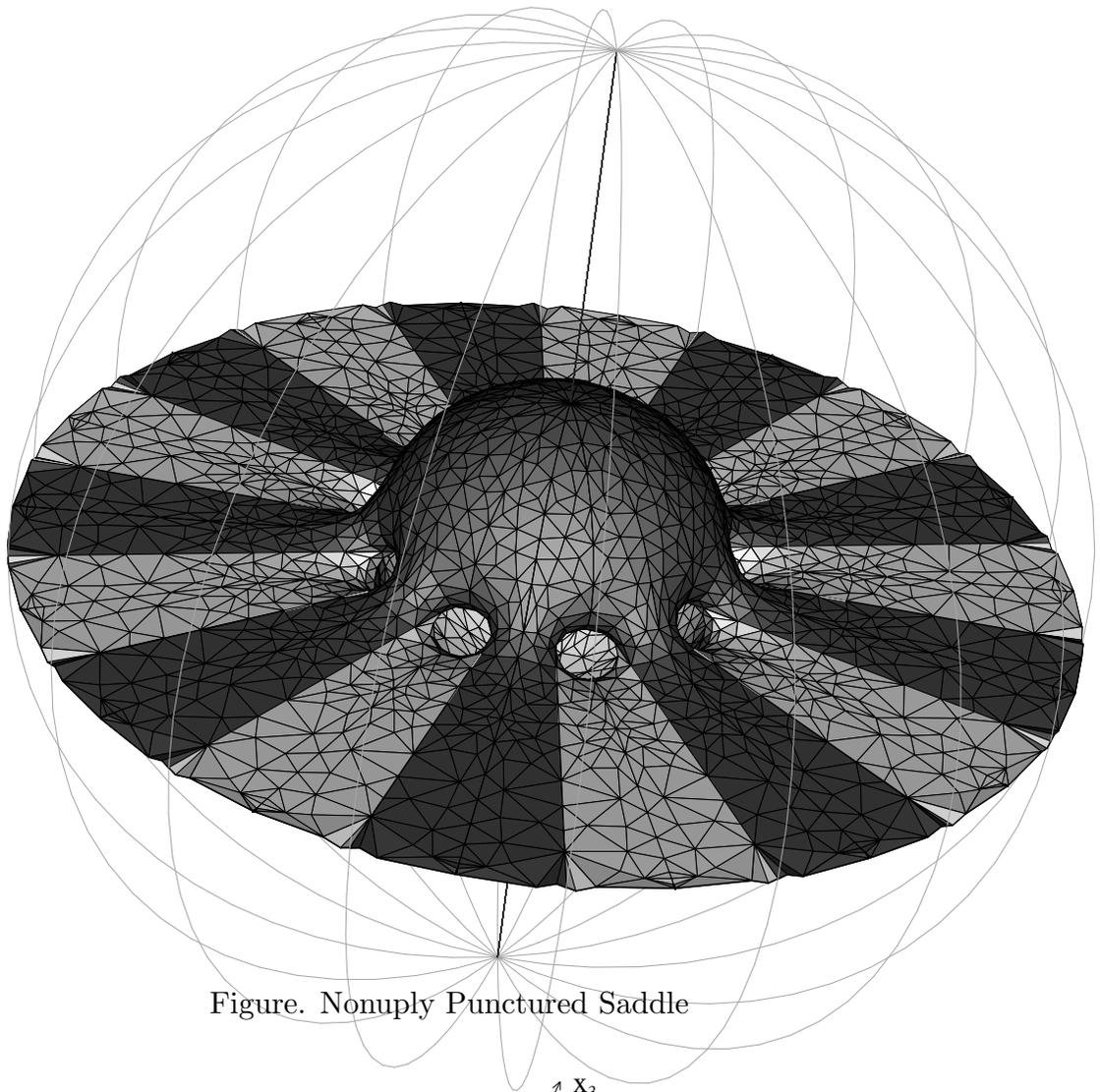

Figure. Nonuply Punctured Saddle

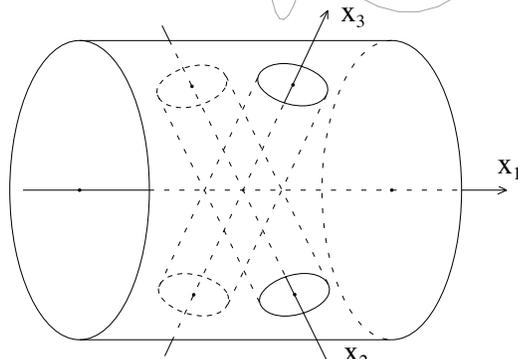

Figure. Initial Surface

To recognize the level-set that develops the interesting singularity, define $R(z, t)$ to be the radius of the radial necks divided by the radius of the longitudinal necks, for the level set $\{u(\cdot, t) = z\}$.

At one extreme, $R(z, t) \to 0$ as the radial necks pinch off. At the other extreme, $R(z, t) \to \infty$ as the longitudinal necks pinch off. At each time $t$, select the level-set $z$ for which $R$ is changing the least, and make it the new zero-set. On this level-set, the two sets of necks are shrinking at around the same rate. Also, rescale the whole picture by a homothety so that the two radii remain roughly constant in size.

This creates a sequence of surfaces which, according to our computation, converges to the following self-shrinking surface $N$. It represents a singularity where all eight necks pinch off simultaneously. Note that $N$ is asymptotic to a double-lobed cone $C$ at infinity, and $N_t \equiv \sqrt{-t} \to N_0 \equiv C$ as $t \to 0^-$.

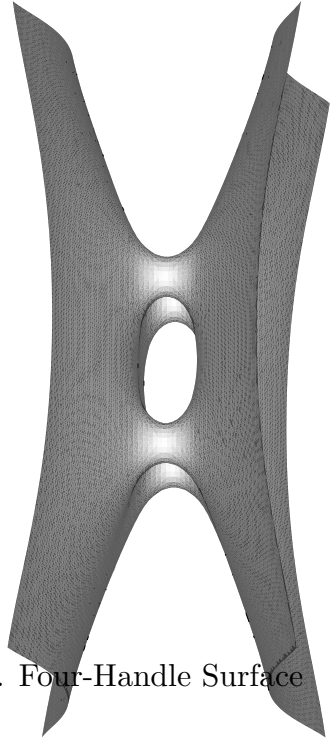

Figure. Four-Handle Surface

Figure. Cone at $t = 0$

Examples with similar topology should exist for any group $G = D_m \times Z_2$, $m \geq 2$. The case $m = 2$ was done in the original paper of Chopp [C].

**Open Problem.** Prove that the $m$-punctured saddle and the $m$-handle surface exist for $m \geq 2$.

## C. Nonuniqueness in $\mathbf{R}^3$

We now attempt to continue the evolution of the four-handle surface past $t = 0$. We obtain

**Result (Angenent-Chopp-Ilmanen).** *The four-handle surface evolves nonuniquely after $t = 0$.*

Here is how we establish this. Observe that the cone $C = N_0$ is the union of two subcones $C^- = C \cap \{x \leq 0\}$ and $C^+ = C \cap \{x \geq 0\}$, each of which is a graph over the $y$-$z$ plane. By standard techniques for the evolution of Lipschitz graphs (see [EH1, EH2]), $C^-$ and $C^+$ have smooth forward evolutions $P_t^-$, $P_t^+$ for $t > 0$, which are smooth graphs over the $y$-$z$-plane. Observe that by the strong maximum principle, $P_t^-$ and $P_t^+$ are disjoint after the first instant. Then the union

$$P_t'' = P_t^+ \cup P_t^-$$

is a forward evolution of the double cone $N_0$ that continues the four-handle evolution. We call $P_t''$ the *two-sheeted* evolution.

Since $P_t^+$ is a graph, it is the unique evolution for its initial condition $N_0$. It follows that $P_t^+$ inherits the symmetry of the equation (MCF), namely

$$P_t^+ = \sqrt{t} \cdot P^+, \qquad t > 0,$$

for some smooth surface $P^+$. (The argument is this: the rescaled flow $\lambda^{-1} \cdot P_{\lambda^2 t}^+$ solves (MCF) with the same initial condition $N_0$, so it must be equal to $P_t^+$.) The same argument works for $P_t^-$.

Therefore $P_t'' = \sqrt{t} \cdot P''$, $t > 0$. We call such an evolution *self-similarly expanding*. (This represents a second solution $a(t)$ of the ODE refered to in section A.)

But there is also a second possible evolution: the two lobes of $C$ can join up near 0 to form a surface of one sheet. Intuitively, this occurs because area decreases sharply by tranching a hole between the cones.

We will construct this one-sheeted evolution by a variational argument which shows that in a certain sense, the two-sheeted evolution is not "minimizing". Inspired by $P_t''$, we impose the self-expander ansatz

$$P_t = \sqrt{t} \cdot P, \qquad t > 0.$$

We also require that $P$ is asymptotic to $C$ at infinity, so that $P_t \to C$ as $t \to 0^+$.

By separation of variables as before, this leads to the elliptic equation

(SE) $$H - \frac{x \cdot \nu}{2} = 0, \qquad x \in P,$$

which says that $P$ expands by homothety under (MCF). This is the Euler-Lagrange equation of the functional

$$\mathbf{K}[P] := \int_P e^{|x|^2/4} \, d\mathcal{H}^{n-1}(x).$$

This functional, in contrast to the self-shrinking functional $\mathbf{J}$, easily yields critical points by minimization; this comes from the convexity of $e^{|x|^2/4}$. In fact we have

**Theorem (Ilmanen [I5]).** *For any closed, homothety invariant set $C$ and any component $E$ of the complement of $C$, there exists a surface $P$ minimizing $\mathbf{K}$ with respect to compact replacements, such that $P$ separates the set*

$$\{x \in E : \operatorname{dist}(x, \partial E) \geq \sqrt{2(n-1)}\}$$

*from the set*

$$\{x \in \mathbf{R}^n \setminus \bar{E} : \operatorname{dist}(x, \partial E) \geq \sqrt{2(n-1)}\}.$$

This theorem is proven by minimizing $\mathbf{K}$ subject to the boundary $Q_R \equiv \partial E \cap \partial B_R$ to obtain a surface $P_R$, then passing a subsequence $P_{R_i} \to P$ as $R_i \to \infty$. This is possible in the space of locally integral $(n-1)$-currents (or locally finite perimeters) using the lower semincontinuity of $\mathbf{K}$ and the compactness theorem of DeGiorgi.

The "separation" condition implies in particular that the limit $P$ is nontrivial, that is, the approximators $P_R$ do not chase off to infinity as $R \to \infty$. This is proven by suitable barriers. It corresponds to the fact that $P_t$ cannot move farther than distance $\sqrt{2(n-1)}$ from $\partial E$ before $t = 1$, because by the maximum principle, $P_t$ cannot collide with the shrinking sphere $\partial B_{R(t)}(x)$ before it disappears at $t = 1$, where $R(t) = \sqrt{2(n-1) - 2(n-1)t}$.

Because of the separation condition, $P_t \to \partial E$ as $t \to 0^+$ in the sense of Hausdorff distance of closed sets, so $P_t$ is a forward evolution of $\partial E$.

Furthermore (by the regularity theory of DeGiorgi, Federer, Almgren, Simons) the singularities of $P$ have codimension at least 7, so in particular $P_t$ is smooth in $\mathbf{R}^3$.

Next we will consider the evolution of the double cone $D_\alpha$ consisting of all points making an angle of $\alpha$ with the $x$-axis, $0 < \alpha < \pi$. We have

**Proposition.** *For each $n \geq 3$, there is a critical angle $\alpha_{crit}(n) \in (0, \pi/2)$ such that for $\alpha \geq \alpha_{crit}$, $D_\alpha$ has both a two-sheeted and a one-sheeted evolution.*

**Proof Sketch.** We have already constructed the 2-sheeted evolution $P_t''(\alpha)$ of $D_\alpha$. Write $P''(\alpha) = P_1''(\alpha)$.

Suppose that $\alpha$ is very close to $\pi/2$. By applying the maximum principle with very large shrinking balls, we can prove that for every $R > 0$, $\varepsilon > 0$, there exists $\alpha(R, \varepsilon) < \pi/2$ such that for all $\alpha \in (\alpha(R, \varepsilon), \pi/2$ we have

$$P''(\alpha) \cap B_R \quad \text{lies within } \varepsilon \text{ of the plane } \{x = 0\}$$

Then by taking $\varepsilon$ very small and $\alpha$ very close to $\pi/2$, we see that the two-sheeted $P''(\alpha)$ cannot be a minimizer of $\mathbf{K}$, because we can replace the two disks $P''(\alpha) \cap B_R$ by a narrow strip in $\partial B_R$ connecting the edge of the disks, thereby greatly reducing $\mathbf{K}$.

This shows $P''(\alpha) \neq P'(\alpha)$, where $P'(\alpha)$ is the minimizer of **K** provided by the Theorem above.

Therefore $D_\alpha$ evolves nonuniquely. By further analysis, we can show that $P'(\alpha)$ really has only one sheet; for details see [ACI]. This completes the proof sketch.

To complete the analysis of the four-handle surface, we need a more quantitative result. Since $D_\alpha$ is rotationally symmetric, we may assume that $P'(\alpha)$ is rotationally symmetric, and reduce (SE) to an ODE. Using Mathematica, we plotted some solutions of this ODE for $\mathbf{R}^3$, and obtained the following graphs. To obtain the self-expanding evolutions, rotate the graphs around the $x$-axis. (Notice that above the critical aperture, there are actually two one-sheeted evolutions of the same topology, but different geometry.)

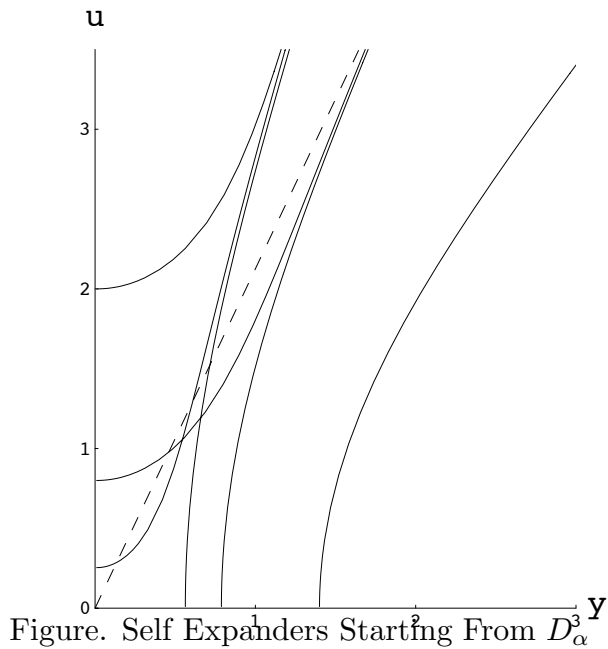

Figure. Self Expanders Starting From $D_\alpha^3$

As a result of this numerical investigation, we found that

$$\alpha_{crit}(3) \approx 66.04°.$$

Returning to the computed four-handle surface, we find by inspection that the cone $N_0$ has "aperture" around $72°$, that is, the cone $D_{72°}$ fits between $N_0$ and the $x$-axis.

Employing the one-sheeted evolution for $D_\alpha$ as a barrier, we can then prove that there exists a one-sheeted evolution

$$P'_t = \sqrt{t} \cdot P', \qquad t > 0,$$

with initial condition $N_0$. The details of this argument are in [ACI]. This establishes (numerically) that the four-handle surface evolves nonuniquely after $t = 0$.

We ask: *which is the preferred evolution?* A good candidate is the minimizer of **K**, since it seems to decrease area fastest, in some weighted sense (this is a suggestion of K. Brakke). But how can we generalize this idea to arbitrary mean curvature flows?

One might ask: what is the probability of various forward evolutions? (Katsoulakis) One way to approach this is to randomize the initial condition or the equation and look at the probability of various outcomes. This seems susceptible to a Monte Carlo simulation.

## D. Level-Set Flow.

Our purpose in sections D-F is to investigate the topological changes that can occur in the flow after the onset of singularities.

To do that, we need another tool, called the *level-set flow*, which is a sort of envelope for all the possible nonunique flows. The level-set flow is due to Chen-Giga-Goto [CGG] and Evans-Spruck [ES1-ES4], based on the idea of viscosity solutions [CIL]. The geometric formulation that I will present in this section appears in Ilmanen [I2], see also Soner [So].

The level-set flow will allow us to distinguish so-called "outermost" flows among all the possible evolutions of a given surface. They seem to be the simplest topologically.

We observe first that smooth flows $M_t$ and $M'_t$ obey the *maximum principle*, that is, if $M_t$ and $M'_t$ are disjoint, they remain disjoint, and in fact the distance between them is nondecreasing. Reason: at the point of closest approach, the curvature vectors cause the surfaces to pull slightly further apart.

We will now make this into a definition of a weak solution, by using smooth flows as test surfaces for a nonsmooth flow. The basic idea is this : the level-set flow is the largest flowing set that does not violate the maximum principle with smooth test flows.

**Definition.** *A family $\{\Gamma_t\}_{t\geq 0}$ of closed sets is a "subsolution" of the mean curvature flow provided that for every compact, smooth mean curvature flow $\{Q_t\}_{t\in[a,b]}$,*

$$if \quad \Gamma_a \cap Q_a = \emptyset, \quad then \quad \Gamma_t \cap Q_t = \emptyset \ for \ a \leq t \leq b.$$

It follows immediately, by using the translation invariance of (MCF), that $\text{dist}(\Gamma_t, Q_t)$ is nondecreasing.

**Proposition.** *For any closed set $\Gamma_0$ there exists a unique* maximal *subsolution $\{\Gamma_t\}_{t\geq 0}$ with initial condition $\Gamma_0$, that is,*

*(a) $\{\Gamma_t\}_{t\geq 0}$ is a subsolution, and*

*(b) If $\{\Delta_t\}_{t\geq 0}$ is any subsolution with $\Delta_0 = \Gamma_0$, then $\Delta_t \subseteq \Gamma_t$ for all $t \geq 0$.*

This is proven by taking the closure of the union of all subsolutions. The maximal subsolution is called the *level-set flow.* It agrees with the definition given in [ES1] and [CGG]. We should remark that the level-set flow has the semigroup property.

The definition would be rather empty without verifying at least some properties in common with the classical mean curvature flow. Here is a basic theorem (for proof, see [I2] and [ES1]).

**Theorem.** *(a) The level-set flow is equal to the smooth flow for as long as the latter exists.*

*(b) (Maximum Principle) Two level-set flows disjoint at time $t = a$ remain disjoint for all times $t \geq a$. (True by definition if one of them is smooth.)*

*(c) (Filling Property) Let $u(\cdot, 0)$ be a Lipschitz function. Then there is a Lipschitz function $u(\cdot, \cdot)$ such that for each $z$, the set $\Gamma_t^z \equiv \{u(\cdot, t) = z\}$ is moving by the level-set flow.*

Property (c) explains the term "level-set flow", and establishes the relationship with the work of Evans-Spruck and Chan-Giga-Goto. The main point of (c) is that no gaps develop between evolving sets that "foliate" space.

By definition, the level-set flow contains all possible mean curvature flows, so: *Nonuniqueness of the mean curvature flow corresponds to "fattening" of the level-set flow.*

For example, the flows $P_t'$ and $P_t''$ for the four-handle surface are both subsolutions, so they lie within the corresponding level-set flow $\Gamma_t$.

By further analysis [ACI], we establish that for $t > 0$, $\Gamma_t$ consists precisely of the region bounded by $P_t'$ and $P_t''$. ($\Gamma_t$ is equal to $N_t$ the four handle surface for $t \leq 0$.)

More generally, let $\Gamma_0$ be the zero-set of a Lipschitz function $u(\cdot, 0)$, and let $u(x, t)$ be the evolution provided by the "filling" property. Then the two sets

$$M_t^+ = \partial\{u(\cdot, t) > 0\}, \qquad M_t^- = \partial\{u(\cdot, t) < 0\},$$

are called "outermost" flows.[2]

---

[2] In what sense do they satisfy the equation? Using the existence and compactness theorems of Ilmanen [I3], one can prove that each of $M_t^\pm$ is the support of a Brakke flow.

In the case of fattening, the outermost flows are distinct; the two sides of the surface peel apart after the onset of fattening.

In the next section, we will see that outermost flows have some special topological properties.

Figure. Level-Set Flow of $D_\alpha$

## E. Topological Changes Past a Singularity

As background for the next section, we will say something about the possible topological changes that can arise during the flow. The results in this section are due to B. White.

Let $\{\Gamma_t\}_{t \geq 0}$ be a level-set flow, and let $E_t = \mathbf{R}^n \setminus \Gamma_t$. For $a \leq b$, define

$$E[a, b] = \cup_{a \leq t \leq b} E_t \times \{t\}.$$

The set $E[a, b]$ is open in $\mathbf{R}^n \times [a, b]$. (Reason: $E[a, b] = \{u \neq 0\} \cap \mathbf{R}^n \times [a, b]$.)

**Theorem (White [W2]).** *The inclusions $E_a \subseteq E[a, b]$, $E[a, b] \supseteq E_b$ induce the following maps, which are injective, surjective and isomorphic as indicated.*

*(a) $\pi_0(E_a) \cong \pi_0(E[a, b])$,*

*(b) $\pi_1(E_a) \twoheadrightarrow \pi_1(E[a, b])$,*

*(c) $H_{n-2}(E[a, b]) \hookleftarrow H_{n-2}(E_b)$,*

*(d) $H_{n-1}(E[a, b]) \cong H_{n-1}(E_b)$.*

The proof of this theorem uses the maximum principle very cleverly, together with technical tools from [I3]. See [W4] for the proof.

If $n = 3$, then combining (b) and (c) we obtain

$$H_1(E_a) \twoheadrightarrow H_1(E[a, b]) \hookleftarrow H_1(E_b), \qquad a \leq b,$$

so we have

**Corollary (White).** *In $\mathbf{R}^3$, $\dim H_1(E_t)$ is nonincreasing.*

Now let $M_t^+$ be a compact outermost flow in $\mathbf{R}^3$. Then $M_t^+$ divides $\mathbf{R}^3$ into two open regions, $U_t^+ = \{u(\cdot, t) > 0\}$ and $U_t^- = \mathbf{R}^3 \setminus U_t^+ \setminus M_t^+$. If $M_t^+$ happens to be smooth, then by the Mayer-Vietoris sequence and Poincare duality,

$$2\operatorname{genus}(M_t) = \dim H_1(M_t^+) = 2 \dim H_1(U_t^+),$$

which is nonincreasing since $\{\mathbf{R}^3 \setminus U_t^+\}_{t \geq 0}$ is a level-set flow, as can be seen by applying the definition directly (or see [I1] for a proof). This shows:

*The genus of a compact outermost flow is nonincreasing in time, whenever it is defined.*

## F. Evolution of Cones in $\mathbf{R}^3$

Suppose that $C$ is a cone in $\mathbf{R}^3$, and consider its evolution. In this section we ask two questions:

*(i) What is the topology of an outermost flow of $C$?*

*(ii) Do all evolutions of $C$ have to be self-similar?*

**(i) Topology.** Let $\Gamma_t$ be the level-set flow of $C$. Since it is unique, we deduce using the parabolic scaling of (MCF) and of $C$ that $\Gamma_t$ itself is self-expanding:

$$\Gamma_t = \sqrt{t} \cdot \Gamma_1, \qquad t > 0.$$

Then the "outermost" flows will also be self-expanding, and so $\partial \Gamma_1$ will solve the elliptic equation (SE), in some sense. We have the following result.

**Proposition (Ilmanen [I5]).** *The boundary of $\Gamma_t$ is smooth for $t > 0$.*

**Proof Sketch.** We will show that we can approximate $\partial\Gamma_1$ from outside by hypersurfaces that are stationary and stable for $\mathbf{K}$.

Fix $R > 0$, $\varepsilon > 0$. Let $S_\varepsilon = B_R \cap \partial\{x : \operatorname{dist}(x, \Gamma_1) < \varepsilon\}$. Consider the minimization problem

(3) $\qquad \min\{\mathbf{K}[P] : P \in \mathbf{I}_2(\mathbf{R}^3),\ \partial P = \partial S_\varepsilon,\ \operatorname{spt} P \subset\subset \mathbf{R}^3 \setminus \Gamma_1\},$

where $\mathbf{I}_2(\mathbf{R}^3)$ is the space of integral 2-currents in $\mathbf{R}^3$.

Since $\{\Gamma_t\}_{t \geq 0}$ solves (MCF) in the level-set sense and is self-expanding, it follows that $\Gamma_1$ solves (SE) in a "viscosity" sense. Therefore, at least formally, $\Gamma_1$ should be a barrier for equation (SE), and by the strong maximum principle, there should exist a minimizer of (3) that is disjoint from $\Gamma_1$.

In [ISZ, I3, I5], this is made precise, yielding a hypersurface $P_\varepsilon \subseteq \mathbf{R}^3 \setminus \Gamma_1$ with $\partial P_\varepsilon = \partial S_\varepsilon$, which is stationary and stable for $\mathbf{K}$. It is smooth by the estimates of Schoen-Simon [ScS]. Passing $\varepsilon \to 0$, the compactness theorem of [ScS] yields a smooth surface $P_R \subseteq \overline{\mathbf{R}^3 \setminus \Gamma_1}$ with $\partial P_R \subseteq \partial \Gamma_1$.

Now $\{\Gamma_t\}_{t \geq 0}$ is a level-set flow, and $\{\sqrt{t} \cdot P_R\}_{t \geq 0}$ is a smooth mean curvature flow with boundary lying in $\Gamma_t$. It follows that $\{\Gamma_t \cup (\sqrt{t} \cdot P_R)\}_{t \geq 0}$ is a "subsolution" as defined in section D. Then, since $\{\Gamma_t\}_{t \geq 0}$ is maximal,

$$P_R \subseteq \Gamma_1.$$

Then $\partial \Gamma_1 \cap B_r = P_R$ (it requires a short argument to see that every point of $\partial \Gamma_1$ is approached by $P_\varepsilon$). Since $R$ was arbitrary, $\partial \Gamma_1$ is a smooth surface solving (SE), which completes the proof of the proposition.

We now turn to as simple result about the topology of $\partial \Gamma_t$. It is in the same vein as the above theorem of White (see also [W1]); we will be able to give a complete proof here.

Let $M$ be a surface that is a boundary component of a three manifold $K$. We call $M$ *incompressible* in $K$ if every embedded curve in $M$ that bounds an embedded disk in $K$, also bounds an embedded disk in $M$. By the Loop Theorem (Hempel [H]), it is equivalent to require that $\pi_1(M)$ injects into $\pi_1(K)$.

If $M$ is not a closed surface, we define the genus of $M$ to be the smallest genus of a closed surface containing $M$.

**Proposition (Ilmanen-White).** *Let $C \subseteq \mathbf{R}^3$ be a cone, and let $\Gamma_t$ be its level-set flow. For a.e. $R > 0$, each component of $\Gamma_t \cap B_R$ is a 2-disk or a 3-ball, and $\partial \Gamma_t$ has genus zero and is incompressible in $\mathbf{R}^3 \setminus \Gamma_t$.*

The import of the proposition is that the outermost flow of a cone is as simple as possible, topologically speaking. Since cones arise as spatial blowups of a singularity, this helps explain (infinitesimally) why the genus of the outermost flow is nonincreasing.

**Example.** Recall the double cone $D_\alpha$ of section C. If $\alpha$ is small, then topologically, $\Gamma_t$ consists of two disks. If $\alpha$ is large, then $\Gamma_t$ consists of two disks and an annulus. The annulus is incompressible in the complement of $\Gamma_t$.

**Proof.** Let $\gamma$ be a loop in $\Gamma_1 \cap B_R$. By classical results for minimal surfaces, $\gamma$ bounds a smoothly immersed, $\mathbf{K}$-minimizing disk $D$ in $B_R$ (using also the convexity of $B_R$ with respect to $\mathbf{K}$). Note that although it is merely immersed, $\sqrt{t} \cdot D$ is nevertheless a "subsolution", except along its boundary. By the same maximality argument as above, $D \subseteq \Gamma_1 \cap B_R$. Therefore $\pi_1(\Gamma_1 \cap B_R) = 0$.

By the previous proposition, for a.e. $R > 0$, each component $K$ of $\Gamma_1 \cap B_R$ is a 2-manifold with boundary or a topological 3-manifold with boundary. In the first case, $K$ is a 2-disk. In the second case, $K$ is a 3-ball with a finite number of interior balls removed. By maximality again, $K$ is a 3-ball.

2. Thus $\partial(\Gamma_1 \cap B_R)$ consists of spheres, so $\partial \Gamma_1 \cap B_R$ has genus zero and $\partial \Gamma_1$ has genus zero.

3. Now suppose that $\gamma \subseteq \partial \Gamma_1$ bounds an embedded disk $D' \subseteq \overline{\mathbf{R}^3 \setminus \Gamma_1}$. Then we may minimize $\mathbf{K}[D]$ in the space of integral 2-currents subject to

$$\partial D = \partial D', \qquad D \subseteq \overline{\mathbf{R}^3 \setminus \Gamma_1}.$$

Since $\partial \Gamma_1$ is smooth and solves (SE), $D$ is smooth and solves (SE), and either $D \subseteq \partial \Gamma_1$ or $D \cap \partial \Gamma_1 = \partial D \cap \partial D_1$.

By the maximality argument above, the former case holds, and $\partial \Gamma_1$ is incompressible in $\overline{\mathbf{R}^3 \setminus \Gamma_1}$. This completes the proof.

By the nature of these methods, we can draw conclusions only about the outermost flows. What about the flows that lie strictly inside the level-set flow?

**Genus Reduction Conjecture.** *Let $M_t$ be any mean curvature flow in $\mathbf{R}^3$. Then the genus of $M_t$ strictly decreases at the moment of a singularity, unless the singularity is a neckpinch or shrinking sphere.*

This covers genus reduction even for "interior" flows, and also says that, except for the cylinder and the sphere, a singularity consumes a certain amount of topology.

It includes the conjecture that every smooth self-shrinker of genus zero is the plane, cylinder, or sphere.

**(ii) Self-Similarity.** One might think that every evolution of a cone is self-similar. Actually, we expect non-self-similar evolutions of any cone that evolves nonuniquely. Here is a heuristic argument that can be made rigorous in some cases.

Suppose that $P_t$, $P_t'$ are two distinct self-expanding evolutions of $C$. We may assume that either $P_1$ or $P_1'$, say $P_1$, is unstable for the functional $\mathbf{K}$, for otherwise we may construct an unstable critical point of $\mathbf{K}$ with initial condition $C$ using the sweepout method of Pitts [Pi]; see also [ScS]. Then we construct a connecting orbit $\{\tilde{M}_s\}_{s \in \mathbf{R}}$ from $P_1$ to some other critical point, say $P_1''$, under the flow for $\mathbf{K}$, namely

$$\frac{\partial}{\partial s} y = \vec{H}(y) + \frac{y}{2}, \qquad y \in \tilde{M}_s, \quad s \in (-\infty, \infty).$$

Under the change of coordinates $t = e^s$, $x = \sqrt{t}\, y$, $M_t = \sqrt{t} \cdot \tilde{M}_s$, this equation becomes the ordinary mean curvature flow,

$$\frac{\partial}{\partial t} x = \vec{H}(x), \qquad x \in M_t, \quad s \in (-\infty, \infty).$$

Then $M_t$ is close to $P_t$ for $t$ near 0 and close to $P_t'$ for $t$ near $\infty$, and is not self-similar.

This argument can be carried out rigorously for the nonunique examples in Lecture 4. In this way we obtain a huge number of non-self-similar evolutions of cones.

Here is a proposal for a rather spectacular example of the failure of self-similarity, due to B. White. Let $M_0$ consist of the union of two coordinate planes, dividing $\mathbf{R}^3$ into four quadrants labeled I, II, III, IV in order around the axis. Consider a smooth approximator $M_0^\varepsilon$ to $M_0$ constructed as follows. Suppose the singular set of $M_0$ is the $z$-axis. Along the line segments $2k < z < 2k + 1$, make narrow slits that connect I to III, and along the line segments $2k + 1 < z < 2k + 2$, make narrow slits that connect II to IV.

Figure. Starting Surface $M_0^\varepsilon$.

Consider the evolution $M_t^\varepsilon$ of $M_0^\varepsilon$. It can be proven that it converges to the second Scherk surface as $t \to \infty$. (See White [W5] for a similar argument.) This Scherk surface, known as the "Scherk tower", is a minimal surface of the same topology and periodic symmetry as the above.

We speculate: as $\varepsilon \to 0$, $M_t^\varepsilon$ converges to a surface $M_t$ such that $M_t$ is smooth for $t > 0$ and approaches the Scherk surface as $t \to \infty$.

This example leads to many questions: what if we vary the spacing of the holes? For example, what if the holes are invariant under a discrete group of homotheties, rather than translations? Or: if there are a finite number of slits, can we get a self-expander corresponding to a "finite Scherk tower" (Ilmanen-Kusner)? All of these questions are subject to computer investigation.

These examples leave open the possibility that the following conjecture might be true. It is motivated by the analogous result for backward blowups, presented in the next section, which follows from a monotonicity formula.

**Conjecture.** *For any mean curvature flow $M_t$ (including singular) and any point $(x_0, t_0)$, the forward blowups defined by*

$$M_t^\lambda = \lambda^{-1} \cdot M_{t_0 + \lambda^2 t}, \qquad t > 0,$$

*converge subsequentially to self-expanders.*

**Remark.** There is a "forward" monotonicity formula involving $e^{|x|^2/4}$, valid on any surface with self-expanding boundary as follows,

$$\partial M_t = Q_t \qquad \text{where } Q_t = \sqrt{t} \cdot Q, \quad t > 0.$$

However, this formula doesn't seem to be much use for proving the above conjecture, since the integral is infinite on the forward evolution of any cone.

# Lecture 3: Flows in $\mathbf{R}^3$, continued

The outline of this lecture:

  G - Monotonicity Formula

  H - Blowup Theorem in $\mathbf{R}^3$

## G. Monotonicity Formula

We will now discuss the precise meaning of Huisken's principle:

   *Singularity formation is modelled by self-shrinking surfaces.*

That is, if we rescale a singularity parabolically, we obtain (subsequentially) a self-shrinking surface.

First we present the monotonicity formula of Huisken [Hu]. For motivation, we mention that if $u$ solves the ordinary heat equation $u_t = \Delta u$ and $\phi$ is a backwards heat kernel, then

$$\frac{d}{dt} \int_{\mathbf{R}^n} \phi |Du|^2 = \int_{\mathbf{R}^n} -2|\Delta u|^2 \leq 0$$

as the reader may easily calculate.

Let $\rho_{x_0,t_0}(x,t)$ be a backwards heat kernel defined on $\mathbf{R}^n$, but with the scaling appropriate to $\mathbf{R}^{n-1}$, that is,

$$\rho_{x_0,t_0}(x,t) = \frac{1}{(4\pi(t_0-t))^{(n-1)/2}} e^{-|x-x_0|^2/4(t_0-t)}, \qquad t < T, \quad x \in \mathbf{R}^n.$$

**Monotonicity Formula (Huisken).** *For any mean curvature flow $M_t$, $0 \leq t < T$, and any $(x_0, t_0) \in \mathbf{R}^n \times [0, T]$, we have the decreasing quantity*

$$\frac{d}{dt} \int_{M_t} \rho_{x_0,t_0}(x,t) \, \mathcal{H}^{n-1}(x)$$
$$= \int_{M_t} -\rho_{x_0,t_0}(x,t) \left| \vec{H}(x,t) + \frac{(x-x_0) \cdot \nu(x,t)}{2(t_0-t)} \right|^2 d\mathcal{H}^{n-1}(x),$$

*for $0 \leq t < t_0$.*

There are well-known analogues for minimal surfaces (see Simon [S1]), for the harmonic map heat flow (Struwe [St]), the equation $u_t = \Delta u + u^p$ (Giga-Kohn [GK]), and the Yang-Mills heat flow (Price [P]).

**Proof.** Assume that $(x_0, t_0) = (0, 0)$ and write $\rho = \rho_{0,0}(x, t)$. First we use the weighted first variation formula (1) with $X = \vec{H}$, then we use the regular first variation formula (in the "wrong" direction!) with $X = D\rho$ to complete the square. Thus we derive for $t < 0$,

$$
\begin{aligned}
\frac{d}{dt} \int_{M_t} \rho &= \int -\rho H^2 + D\rho \cdot \vec{H} + \rho_t \\
&= \int -\rho^2 + 2D\rho \cdot H - \frac{(D\rho \cdot \nu)^2}{\rho} + \left( \operatorname{div}_M D\rho + \frac{(D\rho \cdot \nu)^2}{\rho} + \rho_t \right) \\
&= \int -\rho \left| \vec{H} - \frac{x \cdot \nu}{2t} \right|^2 + 0
\end{aligned}
$$

because, by a miracle that may be verified by direct calculation, the quantity $\operatorname{div}_M(D\rho) + (D\rho \cdot \nu)^2/\rho + \rho_t$ vanishes identically when we choose $\rho$ to be the backwards Gaussian above. (The reader might be interested in finding other functions $\rho$ with this property.) This proves the monotonicity formula.

Define

$$
\Phi_{(x_0, t_0)}(t) = \int_{M_t} \rho_{x_0, t_0}(x, t) \, \mathcal{H}^{n-1}(x).
$$

Observe that $\Phi$ is a unitless quantity that measures the local area of $M_t$ at time $t$ and scale $\sqrt{t_0 - t}$, smeared out by the Gaussian.

The monotonicity formula yields two consequences: first, $\Phi_{x_0, t_0}(t)$ decreases to a well-defined limit as $t \to t_0$; second, the quantity on the right-hand side must approach zero in some sense. The reader will notice that this quantity is precisely the one that vanishes on self-shrinkers. Let us think how we can take advantage of this.

Assume that for some $C > 0$,

$$(4) \qquad \mathcal{H}^{n-1}(M_0 \cap B_R(x)) \le CR^{n-1}, \qquad R > 0, \quad x \in \mathbf{R}^n.$$

This local area bound is obviously satisfied for compact, smooth initial surfaces and all reasonable noncompact ones.

By the monotonicity formula, it follow that

$$(4') \qquad \mathcal{H}^{n-1}(M_t \cap B_R(x)) \le C, \qquad R > 0, \quad x \in \mathbf{R}^n, \quad t > 0,$$

(with a larger $C$) and

$$(5') \qquad \int_0^{t_0} \int_{M_t} \rho_{x_0,t_0} \left| H - \frac{(x-x_0) \cdot \nu}{2(t_0-t)} \right|^2 d\mathcal{H}^{n-1}(x) \, dt < \infty.$$

Define the *rescaled flow* for any $\lambda > 0$,

$$M_t^\lambda = \lambda^{-1} \cdot (M_{t_0+\lambda^2 t} - x_0), \qquad -t_0/\lambda^2 \leq t < 0.$$

where $-x_0$ represents translation and $\lambda^{-1} \cdot$ represents homothety. By parabolic invariance, $\{M_t^\lambda\}_{t \in [-t_0/\lambda^2, 0)}$ is also a mean curvature flow. The quantities in $(4')$ and $(5')$ are scale invariant, so we obtain

$$(4'') \quad \mathcal{H}^{n-1}(M_t^\lambda \cap B_R(x)) \leq C, \qquad R > 0, \quad x \in \mathbf{R}^n, \quad -t_0/\lambda^2 \leq t < 0,$$

and

$$(5'') \qquad \int_{-100}^0 \int_{M_t^\lambda} \rho \left| H - \frac{x \cdot \nu}{-2t} \right|^2 \to 0 \qquad \text{as } \lambda \to 0.$$

Therefore, in a formal sense, we get $H + x \cdot \nu/2 = 0$ in the limit, and $M_t^\lambda$ converges (at least subsequentially) to a self-shrinker $N_t$. Formally, this justifies the slogan "developing singularities are modeled by self-shrinkers", which we used as motivation for the computer searches.

We call this limit surface, if it exists, a *blowup* of $M_t$ at $(x_0, t_0)$.

**Example.** The blowup of a shrinking convex surface at its unique point of disappearance is a self-shrinking sphere [Hu1].

Blowups always exist in the weak sense of Brakke (moving varifolds); see [I6] or White [W2].

Our aim in the next section is to prove that the blowups exist classically for surfaces in $\mathbf{R}^3$.

## H. Blowup Theorem in $\mathbf{R}^3$.

In this section we will establish that in $\mathbf{R}^3$, a blowup surface (at the first singular time) is smooth. This theorem might help to prove rigorously that the self-shrinkers of Chopp actually exist, and (if pushed a little further) could help in improving the partial regularity theory in $\mathbf{R}^3$.

**$\mathbf{R}^3$ Blowup Theorem.** *In $\mathbf{R}^3$, if the evolving surface is embedded, satisfies (4), and has bounded genus, then the blowups are smooth (though not necessarily the convergence).*

**Remarks.** 1. There is a version of this theorem, with a discrete set of singularities, for immersed surfaces in $\mathbf{R}^3$, and for two dimensional surfaces in higher dimensions.

2. It is not know whether the full blowup sequence converges.

**Sketch of Proof.** The proof mostly involves known techniques in the compactness theory of 2-dimensional surfaces, with a few twists. For general background, see Sacks-Uhlenback [SaU], Schoen-Yau [SY], Choi-Schoen [CS], Gulliver [Gu1, Gu2], and many others; see also Ecker [E].

Let $M_t^\lambda$ be the rescaled surface. Write $M^\lambda = M_{-1}^\lambda$. We will show that a subsequence of $\{M^\lambda\}_{\lambda>0}$ passes to a smooth limit. Let us fix $R > 0$ and prove this in $B_R$. It may later be extended to $\mathbf{R}^3$ by diagonalization.

1. First, let accomplish what we can just using the monotonicity formula and general theorems of geometric measure theory; see [S1] for background.

By (4″) and (5″), there exists a sequence $\lambda_i \to 0$ such that

(6)     $\mathcal{H}^{n-1}(M^{\lambda_i} \cap B_R) \le C, \qquad H + x \cdot \nu/2 \to 0 \qquad \text{in } L^2(M^{\lambda_i}).$

After choosing a subsequence, there exists a Radon measure $\mu$ in $B_R$ satisfying

$$\mu = \lim_{i \to \infty} \mathcal{H}^2 \lfloor M^{\lambda_i} \cap B_R.$$

By Allard's Compactness Theorem [S1] and the bound on $\int H^2$, the limit measure $\mu$ will be *integer 2-rectifiable*, that is

$$d\mu = \theta(x) \, d\mathcal{H}^2 \lfloor X,$$

where $X$ is an $\mathcal{H}^2$-measurable, 2-rectifiable set and $\theta$ is an $\mathcal{H}^2\lfloor X$-integrable, integer valued "multiplicity function". The multiplicity function measures how the limit surfaces pile up. (Equivalently, $\mu$ is the mass measure of an integer rectifiable 2-varifold.)

In particular, $\mu$ has an approximate tangent plane almost everywhere, that is, for $\mu$-a.e. $x$ there exists a blowup measure $\nu$ defined by

$$\nu(A) = \lim_{\sigma \to 0} \sigma^{-2} \mu(x + \sigma \cdot A) \qquad \text{for } A \subseteq \mathbf{R}^2,$$

such that $\nu = k\mathcal{H}^2\lfloor P$ for some 2-plane $P$ and integer $k \geq 1$. We will use $T_x\mu$ to denote the plane $P$ whenever it exists. In particular, it is possible to express divergence-form equations weakly on $\mu$.

Also by Allard's Compactness Theorem, the convergence takes the stronger form

$$\int \phi(x, T_x\mu)\,d\mu(x) = \lim_{i \to \infty} \int_{M^{\lambda_i}} \phi(x, T_x M^{\lambda_i})\,d\mathcal{H}^2(x)$$

for every test function $\phi : \mathbf{R}^3 \times G_2(\mathbf{R}^3) \to \mathbf{R}$, where $G_2(\mathbf{R}^3)$ is the Grassmannian of 2-planes in $\mathbf{R}^3$.

We define $H$ for $\mu$ by taking the first variation formula (2) as a definition, and applying the Radon-Nykodym differentiation theorem. The quantity $\int |H|$ is lower semicontinuous under weak convergence, and using the above displayed formula it is not hard to show that $\int |H + x \cdot \nu|^2$ is lowersemicontinouous as well. Therefore $\mu$ solves (SS) weakly, that is,

$$\int -\operatorname{div}_{T_x\mu} X + \frac{x \cdot (T_x\mu)^\perp \cdot X}{2}\,d\mu = 0$$

for all test vectorfields $X \in C_c^1(B_R, \mathbf{R}^3)$, where $(T_x\mu)^\perp$ represents projection onto the normal space of $T_x\mu$. In particular $H \in L^\infty(\mu)$.

We would then like to invoke Allard's local regularity theorem in the form:

**Allard's Theorem (see [S1]).** *There are constants $c > 0$, $\varepsilon_0 > 0$ depending on $n$ with the following property. Let $\mu$ be an integer rectifiable Radon measure in $B_r$, and assume*

$$|H(x)| \le M \qquad \mu\text{-}a.e.\ x \in B_r,$$

*$r \le c/M$, $0 \in \operatorname{spt}\mu$, and*

$$\mu(B_r) \le (1 + \varepsilon_0)\pi r^2.$$

*Then $\operatorname{spt}\mu \cap B_r/2$ is a $C^{1,\alpha}$ manifold (in fact a graph over a domain in some plane).*

By a bootstrap process, we would then obtain that $\operatorname{spt}\mu \cap B_r/2$ is smooth and (SS) is satisfied classically.

2. Unfortunately, we cannot immediately apply Allard's Theorem because the second hypothesis, which says that $\mu$ has roughly one layer, may not be satisfied.

For minimal surfaces, the archetypical counterexample is the catenoid. It satisfies $H = 0$, but on a large scale, there are roughly two layers. The catenoid is not a $C^{1,\alpha}$ graph over any plane because of the hole, which we can make as small as we wish by scaling.

Figure. Catenoid

To get around this difficulty we use special tricks that work only for two dimensional surfaces. The proof is completed by the two following lemmas. Let $A = (A_{ij})$ denote the second fundamental form.

**Lemma 1.** *The $\int H^2$ estimate yields an $\int |A|^2$ estimate on a smaller ball.*

**Lemma 2.** *If $M_i$ is a sequence of smooth, embedded surfaces in $B_R$ with*

$$\mathcal{H}^2(M_i \cap B_R) \le C, \qquad H = f_i + g_i, \qquad \int_{M_i \cap B_R} |A|^2 \le C,$$

where $f_i$ is bounded in $L^\infty$ and $g_i \to 0$ in $L^2$, then a subsequence converges to a Radon measure $\mu$ in $B_R$ such that $N = \operatorname{spt} \mu$ is a $C^{1,\alpha}$ manifold with $H \in L^\infty$. The convergence may be badly behaved at a finite number of "concentration points" of $|A|^2$.

Lemma 1 is based on the Gauss-Bonnet formula with a cutoff function.

Lemma 2 uses a well known concentrated-compactness argument in the theory of surfaces, together with removal of point singularities. It amounts to the following: a singular point can only occur where topology concentrates; in any ball where this happens, $\int |A|^2$ has a certain minimum value $\varepsilon_0$.

Our implementation runs as follows: away from concentration points of $|A|^2$, we use a Lemma of Simon [S2] to make each component of $M_i \cap B_r(x)$ very close to a plane, then pass to limits and apply Allard's Regularity Theorem separately to the limit of each component.

An alternative is to use the weak $|A|^2$ theory of Hutchinson [Hut1, Hut2]. One could also adapt the classical theory of conformally parametrized surfaces, or the blowup argument of Choi-Schoen [CS].

**Proof of Lemma 1.** We will prove the *localized Gauss-Bonnet estimate*: for any surface $M$ immersed in $\mathbf{R}^3$ and any $\varepsilon > 0$,

$$(7) \qquad (1-\varepsilon) \int\limits_{M \cap B_1} |A|^2 \le \int\limits_{M \cap B_2} H^2 + 8\pi g(M \cap B_2) + \frac{96\pi D}{\varepsilon},$$

where

$$D = \sup_{r \in [1,2]} \frac{\mathcal{H}^2(M \cap B_r)}{\pi r^2}.$$

Scaled to radius $R$, this will establish Lemma 1.

After proving this estimate by a convoluted argument involving curve-shortening on $M$ (!) and tube-counting, the author consulted with R. Gulliver, R. Kusner, T. Toro, B. White, and R. Ye, whose suggestions led to the following, much simpler proof.

The Gauss curvature of $M$ is given by

$$2K = H^2 - |A|^2,$$

since $K = \lambda_1 \lambda_2$, $H = \lambda_1 + \lambda_2$, $|A|^2 = \lambda_1^2 + \lambda_2^2$. Then according to the Gauss-Bonnet formula with boundary, we have for a.e. $r > 0$,

$$\int_{M \cap B_r} |A|^2 = \int_{M \cap B_r} H^2 - 4\pi \chi(M \cap B_r) + 2 \int_{M \cap \partial B_r} \tilde{k},$$

where $\tilde{k}$ is the geodesic curvature of the curve $M \cap B_r$ in $M$. The Euler characteristic is given by

$$\chi(M \cap B_r) = 2c(M \cap B_r) - 2g(M \cap B_r) - h(M \cap B_r),$$

where $c(M \cap B_r)$ is the number of components of $M \cap B_r$, $g(M \cap B_r)$ is the genus of $M \cap B_r$ (after capping off the boundary components by disks), and $h(M \cap B_r)$ is the number of components of $M \cap \partial B_r$. Therefore

$$\int_{M \cap B_r} |A|^2 = \int_{M \cap B_r} H^2 - 8\pi c(M \cap B_r) + 8\pi g(M \cap B_r)$$
$$+ 4\pi h(M \cap B_r) + 2 \int_{M \cap \partial B_r} \tilde{k}.$$

The $c$ term is negative and the $g$ term is controlled by hypothesis. The two remaining terms on the right-hand side are essentially boundary terms. They can be unbounded when transversality is breaking down. We will control them by averaging over $r \in [1, 2]$.

The following lower bound on curvature helps us to control the number of components.

**Lemma 3 (Borsuk [Bo], Milnor [ ]).** *For any curve $\gamma$ in $\mathbf{R}^n$,*

$$\int_\gamma |k| \geq 2\pi h(\gamma),$$

*where $k$ is the geodesic curvature of $\gamma$ in $\mathbf{R}^n$.*

Next we calculate the curvature of $\gamma$.

**Lemma 4.** *If $\gamma$ is the transverse intersection of $M^2$ and $N^{n-1}$, then*

$$|k| \leq \frac{|A_M(\dot{\gamma}, \dot{\gamma})| + |A_N(\dot{\gamma}, \dot{\gamma})|}{\sin \alpha},$$

where $\alpha$ is the angle between $M$, $A_M$ and $A_N$ are the second fundamental forms, and $\dot{\gamma}$ is the unit tangent vector.

We leave the proof of Lemma 4 to the reader.

Applying Lemma 4 we obtain for $\gamma = M \cap \partial B_r$,

$$|k| \leq \frac{|A_M| + 1/r}{\sin \alpha}.$$

We continue the main estimate. Set $\phi = \phi(|x|)$, a cutoff function to be specified below. By the above and Fubini's Theorem,

$$
\begin{aligned}
\int_M \phi |A|^2 \, d\mathcal{H}^2 &= \int_0^1 \int_{M \cap \{\phi > t\}} |A|^2 \, d\mathcal{H}^2 \, dt \\
&= \int_0^1 \bigg( \int_{M \cap \{\phi > t\}} H^2 \, d\mathcal{H}^2 - 8\pi c(t) + 8\pi g(t) + 4\pi h(t) \\
&\qquad\qquad\qquad + 2 \int_{M \cap \{\phi = t\}} |\tilde{k}| \, ds \bigg) dt \\
&\leq \int_M \phi H^2 + 8\pi g(M \cap B_2) + \int_0^1 \int_{M \cap \{\phi = t\}} (2|k| + 2|k|) \, ds \, dt.
\end{aligned}
$$

Here $c(t)$ is the number of components of $M \cap \{\phi > t\}$, etc. In the last line we applied Lemma 3, and used the fact that $|\tilde{k}| \leq |k|$. Let $D_M \phi$ be the tangential derivative $T_x M \cdot D\phi$, and note that $|D_M \phi| = |\sin \alpha||D\phi|$. Using Lemma 4,

$$
\begin{aligned}
\text{(last term)} &= \int_{M \cap B_2} 4|D_M \phi||k| \qquad \text{by the co-area formula} \\
&\leq \int_{M \cap B_2} 4(|D\phi| \sin \alpha) \left( \frac{|A_M| + 1}{\sin \alpha} \right) \\
&= \int_{M \cap B_2} 4|D\phi|(|A_M| + 1) \\
&\leq \int_{M \cap B_2} \frac{4|D\phi|^2}{\varepsilon \phi} + \varepsilon \phi |A_M|^2 + 4|D\phi| \\
&\leq \int_{M \cap B_2} C/\varepsilon + C + \varepsilon \phi |A_M|^2,
\end{aligned}
$$

where we choose $\phi$ to be a cutoff function for $B_1$ inside $B_2$ with $|D\phi|^2/\phi \le C$, $|D\phi| \le C$. Plugging back into the main formula, we get

$$(1-\varepsilon)\int \phi|A|^2 \le \int \phi H^2 + 8\pi g(M \cap B_2) + \frac{C}{\varepsilon}\mathcal{H}^2(M \cap B_2).$$

Select the constants more carefully to get (7). This complete the proof of Lemma 1.

**Proof of Lemma 2.** Write $\mu_i = \mathcal{H}^2\lfloor M_i$. Using (6) and Lemma 1 to bound $\int_{M_1 \cap B_R} |A|^2$, pick a subsequence such that $|A|^2\,d\mu_i$ converges to a Radon measure $\sigma$. Define

$$X = \{x : \sigma(\{x\}) > \varepsilon_2\}$$

for $\varepsilon_2 > 0$ to be chosen. These are called "concentration points". Let $N = \operatorname{spt}\mu$. Then:

(a) $X$ is a finite set.

(b) We will show that $N \setminus X$ is smooth. For each $x \notin X$, there is $r = r(x) > 0$ such that $\int_{M^{\lambda_i} \cap B_r(x)} |A|^2 < \varepsilon_2$ for all sufficiently large $i$. We now invoke Simon's Lemma [S2], which implies the following.

**Lemma.** *For any $C > 0$, there is $\varepsilon_2 = \varepsilon_2(C, n)$ with the following property. Let $M$ be a smooth 2-manifold properly immersed in $B_r \subseteq \mathbf{R}^n$, and suppose*

$$\int_{M \cap B_r} |A|^2 \le \varepsilon \le \varepsilon_2, \qquad \mathcal{H}^2(M \cap B_r) \le Cr^2.$$

*Then for each $y \in M \cap B_{r/4}$, there is a 2-plane $H$ such that*

$$M(y) \subseteq \{x : \operatorname{dist}(x, H) \le \delta_1\}, \qquad \mathcal{H}^2(M(y)) \le (1+\delta_2)\pi(r/2)^2$$

*where $M(y)$ is the component of $M \cap B_{r/2}(y)$ containing $y$, and $\delta_1$, $\delta_2 \to 0$ as $\varepsilon \to 0$.*

The original Lemma is more detailed: $M \cap B_{r'(y)}$ is a topological disk for suitable $r' \in [1/2, 3/4]$, and except for a small set, $M(y)$ is a Lipschitz graph over a domain in $H$, with estimates.

The reader should think of this as a kind of $W^{2,2} \to W^{1,\infty}$ Sobolev embedding theorem that just barely fails. That is,

$$\int |A|^2 \sim \int |D^2 u|^2, \qquad \sup |\nu - \nu_H| \sim \|Du\|_{L^\infty},$$

where $\nu_H$ is the normal to $H$. Now $\sup |\nu - \nu_H|$ just fails to be bounded, so $M(y)$ just fails to be a graph.

Now let $y_i \in M_i$ be any sequence converging to $x$, and let $\tilde{\mu}$ be a subsequential limit of $\mathcal{H}^2 \lfloor M(y_i)$. By previous considerations, $\tilde{\mu}$ is a weak solution of (SS) in $B_{r/2}(x)$. Using the Lemma, Allard's Theorem applies to $\tilde{\mu}$ (after making sure that $r \leq c/M$) and so spt $\tilde{\mu} \cap B_{r/4}(x)$ is smooth. By the monotonicity formula for bounded $|H|$ (see [S1]), $\tilde{\mu}(B_{r/2}(x)) \geq c(r/2)^2$. Therefore the number of components of $M_i$ that approach $x$ is bounded, and we may apply this process simultaneously to them all to conclude that $N$ is a finite union of smooth surfaces near $x$.

Since $M_i$ is embedded, by the strong maximum principle we see that these surfaces are disjoint or equal, so $N$ is smooth near $x$. This shows that $N \setminus X$ is smooth.

(c) Since $|H|$ is bounded on $N \setminus X$, the concentration points are removable singularities by Gulliver [Gu1], [Gu2].

This completes the proof of Lemma 2 and of the Blowup Theorem.

The Blowup Theorem leaves open the possibility that $M^{\lambda_i}$ converges to $\mu$ several layers deep, that is,

$$d\mu = k \, d\mathcal{H}^2 \lfloor N$$

where $N$ is smooth, but $k > 1$.

Even worse, near the concentration points there might be small connections between the layers: little catenoid-like necks, or more complicated structures.

Figure. Bad Convergence

**Conjecture.** *$M^{\lambda_i}$ converges smoothly to $N$, with one layer.*

The main difficulty is that the concentration points might skittle around on the surface of $N$ in an uncontrolled way.

If we could prove this conjecture, it would follow that

**Corollary.** *Every end of a self-shrinker is asymptotic either to a smooth cone or to a self-shrinking cylinder.*

To be more precise: let $N$ be a self-shrinker, then there exists $R > 0$ such that $N \setminus B_R$ decomposes into a finite number of ends $U_j$, such that for each $j$, either

(a) As $\lambda \to 0$, $\lambda \cdot U_j$ converges locally smoothly to a cone $C_j$ such that $C_j \setminus \{0\}$ is smooth.

(b) There is a vector $v_j$ such that $U_j - \tau v_j$ converges to the cylinder $\{x : \mathrm{dist}(x, \mathrm{span}\,(v_j)) = \sqrt{2}\}$ as $\tau \to \infty$.

**Remark.** It follows that $N_t = \sqrt{-t} \cdot N$ converges in the sense of Hausdorff distance to $N_0 = (\cup C_j) \cup (\cup\{\tau v_j : \tau \geq 0\}$.

The proof of this Corollary is beyond the present notes. The basic idea is to apply the Blowup Theorem (in its conjectured strong form) to the self-shrinker $N_t$ at a point $x = x_0$, $t = 0$.

**Conjecture (Rigidity of Cylinder).** *If one end of $N$ is asymptotic to a cylinder, then $N$ is isometric to a cylinder.*[3]

---

[3]Angenent has proposed a likely counterexample: the union of a plane and a cylinder, desingularized along the circle of intersection by tiny holes modelled on the Scherck surface.

# Lecture 4:  Nonuniqueness
## in Geometric Heat Flows

We wish to present some rigorous examples of nonuniqueness in higher dimensions,[4] namely $\mathbf{R}^4$ through $\mathbf{R}^7$. This work is due to Angenent-Velazquez-Ilmanen, to appear in [AIV]. The lecture has five parts:

> A - Overview
>
> B - Cones
>
> C - Forward Time
>
> D - Past Time
>
> E - Reduction of Complexity

We would like to thank K. Brakke, D. Chopp, R. Kusner, R. Mazzeo, M. Paolini, R. Schoen, B. White, and others for valuable discussions.

**A. Overview.**  Remarkably enough, the nonuniqueness phenomenon we are about to describe holds almost identically for several other so-called *geometric heat flows*, namely the harmonic map heat flow, the equation $u_t = \Delta u + u^p$, and (under investigation) the Yang-Mills heat flow. Let us present these.

1. The *harmonic map heat flow* is the equation

$$\frac{\partial u}{\partial t} = \Delta u + |Du|^2 u, \qquad x \in \mathbf{R}^n, \quad t \in \mathbf{R},$$

where $u$ takes values in the unit sphere $S^n = \{|u| = 1\}$ in $\mathbf{R}^{n+1}$. It is the gradient flow for the Dirichlet energy $\int |Du|^2\, dx$ subject to this constraint. See [ESa], [Ch], [CS], [St].

2. The semilinear heat equation

$$\frac{\partial u}{\partial t} = \Delta u + u^p, \qquad x \in \mathbf{R}^n, \quad t \in \mathbf{R},$$

which is the gradient flow for the energy $\int |Du|^2/2 - u^{p+1}/(p+1)$. See [GK], etc.

---

[4]Of course, any example may be crossed with $\mathbf{R}^k$.

3. The *Yang-Mills heat flow* is the parabolic equation

$$\frac{\partial}{\partial t}\nabla = -2\operatorname{div} F_\nabla, \qquad x \in \mathbf{R}^n, \quad t \geq 0.$$

Here $\nabla$ is an affine connection, $F_\nabla = (F^\alpha_{ij\beta}) = -\nabla_i\nabla_j + \nabla_j\nabla_i$ is its curvature, and $(\operatorname{div} F)^\alpha_{i\beta} = \nabla_j F^\alpha_{ij\beta}$. This is the gradient flow for the Yang-Mills functional energy $\int |F_\nabla|^2/2$. See [JT].

**Ricci Flow.** These equations (except the $u^p$ equation) were grouped together by Hamilton in his 1986 ICM lecture [H1] because of their geometric character and many deep analogies between them. This class also includes the Ricci flow

$$\frac{d}{dt}g_{ij} = -2R_{ij}, \qquad x \in M, \quad t \geq 0,$$

where $g_{ij}$ is an evolving Riemannian metric on the manifold $M$ and $R_{ij}$ is its Ricci curvature. The Ricci flow has been the subject of intensive study by Hamilton [H2-H7] and others. A natural question is: *can there be nonuniqueness for the Ricci flow, after the onset of singularities?* This question does not yield to the method of analysis presented in this lcture. We will attempt to resolve these difficulties in future work.

The method (in modified form) also applies to nonlinear wave equations, such as the wave map equation $u_{tt} + |u_t|^2 = \Delta u + |Du|^2 u$; see Shatah-Tahvildar-Zadeh [ST].

**Main Results.** We start with a simple example, which appears in the back of Brakke's book [B]. It is a network with triple junctions moving by mean curvature in the plane, which shrinks self-similarly to a cross, and then expands self-similiarly – in more than one way. The angles at the junctions are maintained at $120°$ except when junctions merge.

Figure. Nonuniqueness in the Plane

However, we desire an example which is smooth, except at the point $x = 0$, $t = 0$. Inspired by Brakke's example and the four-handle example of Lecture 2, we impose the following ansatz, which says that the flow is invariant under parabolic rescaling:

(a) *For $t < 0$, the solution is self-shrinking.*

(b) *At $t = 0$, the solution becomes a cone $C$.*

(c) *For $t > 0$, the solution is self-expanding.*

To enable a rigorous analysis, we also impose:

(d) *The solution is invariant under rotation by the group $SO(p) \times SO(q)$ acting on $\mathbf{R}^n = \mathbf{R}^p \times \mathbf{R}^q$, $p$, $q \geq 2$.*

We have expressed this ansatz for (MCF); for the other equations, "cone" should be replaced by "homothety invariant map" and $SO(p) \times SO(q)$ should be replaced by $SO(n)$.

This ansatz has the effect of reducing (MCF) to an ODE that can be analyzed. Our principal result is the following.

**Theorem [AIV].** *The following equations have solutions satisfying (a)-(d), which evolve nonuniquely after an isolated singularity at $x = 0$, $t = 0$.*

(i) *The mean curvature flow, for hypersurfaces in $\mathbf{R}^4$ through $\mathbf{R}^7$.*

(ii) *The harmonic map heat flow, with domain $\mathbf{R}^3$ through $\mathbf{R}^6$.*

(iii) *The equation $u_t = \Delta u + u^p$ in $R^n$, in the supercritical range*

$$n \geq 3, \qquad p_{crit} < p < p_+,$$

*where $p_{crit} = (n+2)/(n-2)$ is critical for the Sobolev inequality, and $p_+$ is given by*

$$p_+ = \begin{cases} \infty & n \leq 10 \\ \dfrac{n^2 - 8n + 4 + 8\sqrt{n-1}}{(n-10)(n-2)} & n \geq 11 \end{cases}$$

**Remarks.** 1. The explanation of these mysterious ranges is the existence of a finite energy, unstable, static "cone" (i.e. homothety invariant solution). We shall say more about this below.

2. Nonuniqueness for motion by curvature with a time-varying driving term was discovered by Soner-Souganidis [SS]. This result was refined to the equation $\partial x/\partial t = \vec{H} + \nu$ by Bellettini-Paolini [BP].

3. Troy [T], Budd-Qi [BQ] and Lepin [L] studied the self-shrinking solutions for the $u^p$ equation.

4. Coron [Co] found nonuniqueness for harmonic map heat flow (with singular initial data). Our examples satisfy the weak form of the harmonic map heat flow not only with respect to the usual range variations, but also with respect to domain variations. In particular they satisfy the monotonicity formula and usual energy inequality. See Chen-Li-Lin [CLL], Feldman [F], Freire [Fr]. Assuming these latter inequalities, Freire shows uniqueness in dimension $n = 2$.

For simplicity, we will restrict ourselves to discussing the mean curvature flow. We now start a gradual buildup to this theorem.

**B. Cones.** For each $p, q \geq 2$ with $p + q = n$, there is a rotationally invariant stationary hypercone $C_{p,q}$ defined by

$$C_{p,q} = \{x = (y, z) \in \mathbf{R}^p \times \mathbf{R}^q = \mathbf{R}^n : |y|^2/(p-1) = |z|^2/(q-1)\}.$$

The reader may verify by direct calculation that these cones satisfy $H = 0$. (In the case $p = q$ this is clear by symmetry.)

It is well known (see Simons [Si], Bombieri-DeGiorgi-Giusti [BDG], Lawlor [L]), that

- $4 \leq n \leq 7 : C_{p,q}$ is unstable for the area functional (with respect to smooth perturbations), and therefore non-minimizing,

- $\{p, q\} = \{6, 2\} : C_{p,q}$ is stable, and minimizing on one side, but non-minimizing on the other,

- $n \geq 8, \{p, q\} \neq \{6, 2\} : C_{p,q}$ is a minimizer of area (with respect to compact replacements).

**Remarks.** 1. In the harmonic map case, the "cone" is the map $x/|x|$. In the $u^p$ case, the "cone" is the function $A_{p,n}|x|^{-2/(p-1)}$ for $p > n/(n-2)$, where $A_{p,n}$ is a constant.

2. There is also an example of nonuniqueness involving the cone $C_{6,2}$. The mechanism is quite different, however. This was the first nonuniqueness example discovered, but unfortunately we will not have time to cover it.

3. There is no stationary cone in $\mathbf{R}^3$ (other than the plane), even though $n = 3$ lies within the range of dimensions $5 - 2\sqrt{2} < n < 5 + 2\sqrt{2}$ for which such a cone would be unstable. As a result, in $\mathbf{R}^3$, the method presented here produces no embedded, rotationally symmetric examples of nonuniqueness, but there are likely to be immersed, rotationally symmetric examples (in addition to the non-rotationally symmetric example in Lecture 2).

As a warmup to our nonuniqueness analysis, we present the following general fact.

**Proposition (Ilmanen, [I6]).** *Let $C$ be a stationary cone of the form $C = \partial E$, where $E$ is an open set of locally finite perimeter. If $C$ is not area-minimizing, then $C$ evolves nonuniquely.*

**Proof Sketch.** Since $C$ is stationary, one evolution is given by $M_t \equiv C$, $t \geq 0$. On the other hand, by the Theorem quoted in Lecture 2, section C, there exists another evolution of the form $P_t = \sqrt{t} \cdot P$, where $P$ is a minimizer of the functional $\mathbf{K}[P] = \int_P e^{|x|^2/4}$. If $P = C$, then by blowing up at the origin, $\mathbf{K}$ converges to the area functional, and we find that $C$ minimizes area. This establishes the Proposition.

In section C we will present a more precise statement: *$C_{p,q}$ has an infinite number of forward evolutions.*

We give a heuristic explanation of this. The stationary cone $C_{p,q}$ is unstable, which means that there exists a smooth, compactly supported perturbation $X$ of $C_{p,q}$ such that the second variation of area with respect to $X$ is negative (see the figure). Since (MCF) is the gradient flow for area, such a perturbation will grow under the flow.

Figure. A Totally Unstable Cone

We may arrange that $X$ is supported away from the origin. Then we scale $X$ by a homothety until it is disjoint from its original position, and reverse its sign across $C_{p,q}$. (That is, $X_{new}(x) = -\lambda X(x/\lambda)$, where we assume that $X$ is perpendicular to $C_{p,q}$.) This creates a tendency for the cone to "break" on the opposite side. The relative effect of the new perturbation is the same as the original one, except speeded up by a factor $\lambda^{-2}$. By iterating this procedure, we may obtain conflicting tendencies of ever greater speed all the way down.

As we shall see in section C, the cone can indeed "break" on either side with any number of kinks.

**A Nearby Minimal Surface.** Before looking at dynamics, we will describe a smooth minimal surface near $C_{p,q}$ which exhibits infinitely many kinks.

We pause to set up some notation. Suppose that a surface $M$ is invariant under the natural action of $SO(p) \times SO(q)$ on $\mathbf{R}^n = \mathbf{R}^p \times \mathbf{R}^q$. Write $x = (y, z)$ respecting this decomposition. Then $M$ has the form

$$M = M(\gamma) \equiv \{(y, z) : (|y|, |z| \in \gamma\},$$

where $\gamma$ is a curve in $\mathbf{R}^2$. Write $(r, u) = (|y|, |z|)$, and by abuse of notation, write $x = (r, u)$. If $\gamma$ is the graph of a function $u = u(r)$, write $M = M(u)$.

Let $C_\alpha$ be the cone $M(\{u = r \tan \alpha\})$. Define

$$\lambda_s = \tan \alpha_s = \sqrt{\frac{q-1)}{p-1}},$$

so $C_{p,q} = C_{\alpha_s}$.

We now prepare to reduce PDEs to ODEs via the rotational symmetry. The principal curvatures of $M$ with respect to a unit normal $\nu$ are

$$k, \underbrace{\frac{-\cos\theta}{u}, \dots, \frac{-\cos\theta}{u}}_{q-1}, \underbrace{\frac{\sin\theta}{r}, \dots, \frac{\sin\theta}{r}}_{p-1}$$

where $\theta$ is defined by $\dot\gamma = (\cos\theta, \sin\theta)$, $\nu = (-\sin\theta, \cos\theta)$, and $k$ is the curvature of $\gamma$ with respect to $\nu$. Then we obtain:

$$(8) \qquad H = k + \frac{(q-1)\sin\theta}{u} - \frac{(p-1)\cos\theta}{u}, \qquad x \in M.$$

Therefore the minimal surface equation $H = 0$ becomes the following ODE for $\gamma$,

$$k + \frac{(q-1)\sin\theta}{u} - \frac{(p-1)\cos\theta}{u} = 0, \qquad x \in \gamma.$$

Except at isolated vertical points, $\gamma$ is locally the graph of a function $u = u(r)$ solving the following equivalent ODE, using $(\cos\theta, \sin\theta) = \pm(1, Du)/\sqrt{1 + |Du|^2}$.

$$(9) \qquad \frac{u_{rr}}{1 + u_r^2} + \frac{p-1}{r} u_r - \frac{q-1}{u} = 0, \qquad r \geq 0.$$

Now we turn to a smooth minimal surface that is a "companion" to $C_{p,q}$. An analogous example was found by Schoen-Uhlenbeck [ScU2] for harmonic maps.

**Proposition.** *For $4 \leq n \leq 7$, there exists a smooth, $SO(p) \times SO(q)$-invariant minimal surface $Q$ which is asymptotic to $C_{p,q}$ at infinity, and crosses it infinitely many times.*

In fact, $Q = M(v)$, where $v = v(r)$ is a function defined for all $r \geq 0$, and $(v/r, v_r) \to (\lambda_s, \lambda_s)$ as $r \to \infty$. Note that for any $\lambda > 0$, $\lambda \cdot Q$ is also a minimal surface, so we may assume $v(0) = 1$.

In the picture, the reader will observe that $Q$ crosses $C_{p,q}$ an infinite number of times as we pass to infinity.

Figure. A Wiggly Companion

**Proof.** Begin integrating (9) with the initial condition

$$u(0) = 1.$$

By analyzing the singular point at $r = 0$, we find that there is a unique smooth solution defined for small $r$, and that $u_r(0) = 0$ is forced. Furthermore, $u_r > 0$ for small $r > 0$. We now change coordinates by

$$\eta = \log r, \qquad X = u/r, \qquad Y = u_r,$$

rendering the equation automonous, and yielding

$$X_\eta = Y - X, \qquad Y_\eta = -\frac{(p-1)(1+Y^2)(\lambda_s^2 - XY)}{X}.$$

Examining the phase plane, it is easily seen that the solution remains in the first quadrant, and spirals in to the unique fixed point

$$(X, Y) = (\lambda_s, \lambda_s).$$

The details are left to the reader. This completes the proof.

Figure. Phase Plane

In the next two sections, we will look at the forward and backward evolution of cones. (The backward problem is more fragile.) We will ultimately produce

(1) a continuous family of self-expanding solutions (for $t > 0$),

(2) a discrete sequence of self-shrinking solutions (for $t < 0$),

which we will match up at $t = 0$.

**C. Forward Time.** First we construct various expanders for $t > 0$. We recall the self-expanding equation

(SE) $$H(x) - \frac{x \cdot \nu}{2} = 0, \qquad x \in P,$$

which implies that $P_t \equiv \sqrt{t} \cdot P$, $t > 0$, is a mean curvature flow.

Let $P = M(\gamma)$. Using (8) and noting that

$$\frac{x \cdot \nu}{2} = \frac{(-r\sin\theta, u\cos\theta)}{2},$$

we find that $\gamma$ solves the ODE

$$k + \left(\frac{p-1}{r} + \frac{r}{2}\right)\sin\theta - \left(\frac{u}{2} + \frac{q-1}{u}\right)\cos\theta = 0, \qquad x \in \gamma.$$

This is equivalent to the following equation wherever $\gamma$ can be expressed as a graph $u = u(r)$:

$$(10) \qquad \frac{u_{rr}}{1 + u_r^2} + \left(\frac{p-1}{r} + \frac{r}{2}\right)u_r - \frac{u}{2} - \frac{q-1}{u} = 0, \qquad r > 0.$$

Our next task is to describe the solutions of this ODE. One solution is $u = \lambda_s r$. That is, the static cone $C_{p,q}$ is expanding by homothety, without anyone noticing.

**Lemma.** *For each $a > 0$, there exists a minimal surface $P^a = M(u^a)$, where $u_a$ is the solution of (10) with initial condition*

$$u_a(0) = a.$$

*$u_a$ satisfies $(u_a)_r(0) = 0$, is defined for all $r \geq 0$, and is asymptotic to some ray $u = \lambda(a)r$ where $0 < \lambda(a) < \infty$, in the sense that $(u_a/r, (u_a)_r)$ converges to $(\lambda(a), \lambda(a))$ as $r \to \infty$.*

Let $\gamma_a$ denote graph $(u_a)$. Define the angle $\alpha(a)$ by $\lambda(a) = \tan(\alpha(a))$. The lemma implies that the flow assumes a cone as initial condition:

$$P_t^a \equiv \sqrt{t} \cdot P^a \to C_{\alpha(a)}, \qquad t \to 0^+.$$

**Proof Sketch.** The proof is similar to that of the existence of the wiggly companion $M(v)$. We perform the same change of variables and obtain

$$X_\eta = Y - X, \qquad Y_\eta = -(1 + Y^2)\left(\frac{(p-1)(\lambda_s^2 - XY)}{X} + \frac{e^{2\eta}}{2}(Y - X)\right).$$

As before the condition $0 < X < \infty$, $0 < Y < \infty$ is maintained, and in fact $(X, Y)$ never leaves any square with corners on the locus $X = Y$, $XY = \lambda^2$. The extra terms only help to maintain these conditions. However, there is so much friction that $(X, Y)$ can slowly grind to a halt at some point $(\lambda, \lambda)$ other than $(\lambda_s, \lambda_s)$. The details will appear in [AIV].

Figure. Phase Plane

Next, we will investigate the behavior of $\alpha(a)$ as $a$ varies. We will match asymptotics from two regions:

(i) the region $|x| \ll 1$, where $P^a$ is nearly a minimal surface.

(ii) the region $|\alpha(x) - \alpha_s| \ll 1$, where we linearize (10) over the cone.

*(i) Behavior near the origin.*

Near the origin, we ignore the $x/2$ term in (SE) to obtain the minimal surface equation $H = 0$. We expect that $u_a$ will be well approximated by the wiggly curve $av(r/a)$ in the region $r \ll 1$.

Let us examine how $av(r/a)$ approaches the cone. First we will look at $v$. By analyzing the regular singular point of (9) at $r = \infty$ (or equivalently, studying the fixed point $(\lambda_s, \lambda_s)$), we deduce

$$(11) \qquad v(r) \approx \lambda_s r + r^{-\beta}(A_1 \cos(\mu \log r) + A_2 \sin(\mu \log r)), \qquad r \gg 1,$$

where

$$(12) \qquad \beta = \frac{n-3}{2}, \qquad \mu = \frac{\sqrt{8 - (n-5)^2}}{2},$$

and $A_1$, $A_2$ are constants that are uniquely determined. Then

$$u_a \approx av(r/a) \approx \lambda_s r + a(r/a)^{-\beta}(A_1 \cos(\mu \log(r/a)) + A_2 \sin(\mu \log(r/a)))$$

in the regime $a \ll r \ll 1$.

This analysis is valid precisely in the range $4 \leq n \leq 7$; outside of this range the behavior is not oscillatory. This corresponds to the fact that the cone is unstable for area in these dimensions.

*(ii) Behavior near the Static Cone.*

We linearize equation (10) over $C_{p,q}$. Write $u(r) = \lambda_s r + h(r)$. Then, up to an error that is quadratic in $h$ and its derivatives, $h$ solves the equation

$$(13) \qquad \frac{h_{rr}}{1+\lambda_s^2} + \left(\frac{p-1}{r} + \frac{r}{2}\right) h_r + \left(-\frac{1}{2} + \frac{p-1}{r^2}\right) h = 0.$$

Figure. Solutions of Linearized Equation

We find that $h$ has the asymptotics given by (11) as $r$ approaches the regular singular point $r = 0$. This is not surprising, since we have completed a "commutative diagram" of approximations. Therefore we can select a basis $h_1$, $h_2$ for the solutions of (13) by imposing the normalization

$$h_1 \approx r^{-\beta} \cos(\mu \log r), \qquad h_2 \approx r^{-\beta} \sin(\mu \log r), \qquad r \ll 1.$$

At the irregular singular point $r = \infty$, the solution $h$ has the asymptotics

$$h \sim C_1' r + C_2' e^{-(1+\lambda_s^2)r^2/4}/r.$$

The decaying exponential stabilizes the solution, and explains why every initial condition is asymptotic to a cone at infinity.

Therefore there exist constants $\lambda_1$, $\lambda_2$ such that

$$(14) \qquad h_1 \approx \lambda_1 r, \qquad h_2 \approx \lambda_2 r, \qquad r \gg 1.$$

*(iii) Matching*

Now we scale $v$ and $h$ to match in the intersection of regions (i) and (ii). Let $a$, $C_1$, $C_2$ be undetermined constants and set

$$av(r/a) \approx \lambda_s r + C_1 h_1(r) + C_2 h_2(r), \qquad a \ll r \ll 1.$$

This becomes upon substitution,

$$\lambda_s r + a(r/a)^{-\beta}(A_1 \cos(\mu \log(r/a)) + A_2 \sin(\mu \log(r/a)))$$
$$\approx \lambda_s r + C_1 r^{-\beta} \cos(\mu \log r) + C_2 r^{-\beta} \sin(\mu \log r),$$

which yields

$$C_1 \approx a^{\beta+1}(A_1 \cos(\mu \log a) - A_2 \sin(\mu \log a)),$$
$$C_2 \approx a^{\beta+1}(A_2 \cos(\mu \log a) + A_1 \sin(\mu \log a)),$$

and using (14), we estimate the limiting slope

$$
\begin{aligned}
(15) \quad \lambda(a) &\approx \lambda_s + C_1\lambda_1 + C_2\lambda_2 \\
&\approx \lambda_s + a^{\beta+1}((\lambda_1 A_1 + \lambda_2 A_2)\cos(\mu \log a) \\
&\qquad\qquad + (-\lambda_1 A_2 + \lambda_2 A_1)\sin(\mu \log a)) \\
&\equiv \lambda_s + a^{\beta+1}(D_1 \cos(\mu \log a) + D_2 \sin(\mu \log a)).
\end{aligned}
$$

Figure. Matching

Summarizing the above, we have

**Lemma 1 (Self-Expanders).** *For each $a > 0$, the self-expander $P^a$ has the asymptotics given in (i) and (ii). The limiting angle $\alpha(a)$ is given asymptotically by*

$$\tan \alpha(a) \approx \lambda_s + a^{(n-1)/2}(D_1 \cos(\mu \log a) + D_2 \sin(\mu \log a)), \qquad a \ll 1,$$

*where $\mu = \sqrt{8 - (n-5)^2}/2$. Also $\alpha(a) \to \pi/2$ as $a \to \infty$.*

The proof will appear in [AIV].

Figure. Graph of $\alpha(a)$

By counting the set $\{\alpha(a) = \alpha\}$, we obtain the following.

**Corollary.** *The cone $C_{p,q}$ has infinitely many self-expanding forward evolutions. The cone $C_\alpha$ has a very large $L(\alpha)$ number of forward evolutions, given by*

$$L(\alpha) = -\frac{2\mu}{(n-1)\pi} \log(\tan \alpha - \lambda_s) + O(1).$$

**D. Past Time.** Next, we try to concoct a self-shrinker for $t < 0$ that matches up with a cone at $t = 0$. We recall the self-shrinking equation

$$\text{(SS)} \qquad\qquad H(x) + \frac{x \cdot \nu}{2} = 0, \qquad x \in N.$$

If $N = M(\delta)$, then the curve $\delta$ solves the ODE

$$k + \left( \frac{p-1}{r} - \frac{r}{2} \right) \sin \theta + \left( \frac{u}{2} - \frac{q-1}{u} \right) \cos \theta = 0, \qquad x \in \delta,$$

which is equivalent to the following ODE wherever $\gamma$ is locally a graph,

$$\text{(16)} \qquad \frac{u_{rr}}{1 + u_r^2} + \left( \frac{p-1}{r} - \frac{r}{2} \right) u_r + \frac{u}{2} - \frac{q-1}{u} = 0, \qquad r > 0.$$

We must describe the solutions of this ODE. Define $\delta_a$ to be the curve solving (16) with initial condition

$$u(0) = a, \qquad a > 0.$$

In constrast to the self-expanding case, the solutions are wildly unstable for most values of $a$, and will not approach a cone at infinity (the transport terms have the wrong sign). To understand this behavior, we embark again on matched asymptotics in the regions

    (i) near the origin

    (ii) near the cone.

(i) As before, we expect the solution $\delta_a$ of (16) to behave like the wiggly minimal surface $M(av(r/a))$ as long as it stays in the region $|x| \ll 1$. In particular, for small $a$, $\delta_a$ will approach closely to the cone in the ocillatory manner given by (11).

(ii) We linearize equation (16) over the cone. Write $u(r) = \lambda_s r + g(r)$. This time we obtain the equation

$$\text{(17)} \qquad \frac{g_{rr}}{1 + \lambda_s^2} + \left( \frac{p-1}{r} - \frac{r}{2} \right) g_r + \left( \frac{1}{2} + \frac{p-1}{r^2} \right) g = 0.$$

Again, the solution $g$ has the asymptotics given by (11) as $r$ approaches 0, so we can define a basis $g_1$, $g_2$ for the solutions of (17) by requiring

$$g_1 \approx r^{-\beta}\cos(\mu\log r), \qquad g_2 \approx r^{-\beta}\sin(\mu\log r), \qquad \text{for } r \ll 1.$$

As $r \to \infty$, the solution $g$ now has the asymptotics

$$g \sim C_1' r + C_2' e^{(1+\lambda_s^2)r^2/4}/r.$$

Figure. Solutions of Linearized Equation

The consequences for the nonlinear equation are the following. If we are fortunate enough that $C_2'$ equals 0, then $\delta_a$ approaches a cone at infinity with slope $\lambda_s + C_1'$. Otherwise, $\delta_a$ very quickly leaves the cone, going either up or down. Once it gets a little away from the cone, the linearization loses its validity. Eventually $\gamma_a$ "turns", that is, the inclination $\theta$ passes $\pi/2$ or 0 and $\delta_a$ ceases to be an increasing graph. Subsequently, $\delta_a$ wanders all over $\mathbf{R}^2$, seemingly at random. Pictures of such curves appear in Angenent [A5].

Accordingly, we define the solution $g_3$ of (17) and constants $B_1$, $B_2$ by requiring that

(18) $$g_3 \approx r \quad \text{for } r \gg 1; \qquad g_3 = B_1 g_1 + B_2 g_2.$$

(iii) Now we do the matching. Let $a$, $C$ be undetermined constants and set

$$av(r/a) \approx \lambda_s r + C g_3(r), \qquad a \ll r \ll 1.$$

This becomes upon substitution,

$$\lambda_s r + a(r/a)^{-\beta}(A_1\cos(\mu\log(r/a)) + A_2\sin(\mu\log(r/a)))$$
$$\approx \lambda_s r + C(B_1 r^{-\beta}\cos(\mu\log r) + B_2 r^{-\beta}\sin(\mu\log r)),$$

which yields

$$CB_1 \approx a^{\beta+1}(A_1\cos(\mu\log a) - A_2\sin(\mu\log a)),$$
$$CB_2 \approx a^{\beta+1}(A_2\cos(\mu\log a) + A_1\sin(\mu\log a)),$$

or equivalently

$$C(B_1 + iB_2) \approx a^{\beta+1} e^{i\mu \log a}(A_1 + iA_2).$$

This has a discrete sequence of solutions given by

$$\mu \log a_k \approx E + \pi k, \qquad C_k \approx (-1)^k D a^{\beta+1}, \qquad k \in \mathbf{Z}, \quad k \gg 1,$$

where $D$, $E$ are defined by

$$B_1 + iB_2 = e^{iE}(A_1 + iA_2) \qquad D = \frac{\sqrt{A_1^2 + A_2^2}}{\sqrt{B_1^2 + B_2^2}}.$$

We now change $E$ by a multiple of $2\pi$ so that $k$ becomes the number of intersections of $\delta_a$ with the ray $u = \lambda_s r$ (for large $k$).

We now interpret this result. As we gradually decrease $a$ toward 0, because of (i) the solution $\delta_a$ will acquire more and more intersections with the cone $u = \lambda_s r$. As $a$ passes $a_k$, the curve $\delta_a$ will acquire a new intersection as it switches from going "up" to going "down" or vice versa (that is, as the sign of $C_2'$ changes). When $a = a_k$, $\delta_a$ goes neither up nor down, and $\delta_a$ is a complete graph asymptotic to a cone at infinity.

Solving for $a_k$ and observing from (18) that $\lambda_k = \lambda_s + C_k$, we have the following result. Here $F$ and $G$ are given by $F = e^E$, $G = DF^{\beta+1}$.

**Lemma 2 (Self-Shrinkers).** *For $k$ sufficiently large, there exist numbers $a_k > 0$, $\alpha_k \in (0, \pi/2)$ with the asymptotic behavior*

$$a_k \approx F\tau^{-k}, \qquad \tan \alpha_k \approx \lambda_s + G\sigma^{-k},$$

*where*

$$\tau = e^{-\pi/\mu}, \qquad \sigma = \tau^{\beta+1} = e^{-\pi(n-1)/2\mu},$$

*and rotationally symmetric self-shrinkers $N^k = M(\delta_{a_k})$, such that*

*(a) The curve $\delta_{a_k}$ is a complete graph that intersects $u = \lambda_s r$ in $k$ distinct points.*

*(b) $N^k$ is asymptotic to the cone $C_{\alpha_k}$ at infinity.*

The proof will appear in [AVI]. It is based on the asymptotics and shooting argument described above.

Figure. Typical Self-Shrinkers

It follows from (b) that

$$N_t^k \equiv \sqrt{-t} \cdot N^k \to C_{\alpha_k} \qquad \text{as } t \to 0^-.$$

Then we can conveniently look up the number of forward evolutions of $C_{\alpha_k}$ in Lemma 1 of the previous section. We substitute the value of $\alpha_k$, and the various logarithms cancel to yield the following.

**Corollary (Nonuniqueness)** *For large $k$, the self-shrinker $N_t^k$ with $k$ intersections has $k + O(1)$ different continuations of the form $P_t^a$.*

**Remark.** We have only considered continuations that are given by graphs $u = u(r)$. There is another family of continuations given by graphs $r = r(u)$.

Figure. Nonunique Evolution

### E. Reduction of Complexity.

**Proposition.** *The number of intersections with $C_{p,q}$ drops by at least one after the singularity.*

We prove this via the intersection theory of Angenent [A1], which says: if two evolving curves solve the same (smooth, nonlinear) parabolic equation, then

(a) the number of intersections of the two curves is finite after an instant,

(b) the number of intersections is nonincreasing, and strictly decreases whenever the curves are nontransverse.

Note that if $M_t = M(\gamma_t)$ is a smooth mean curvature flow, then $\gamma_t$ solves

$$(19) \qquad \frac{\partial}{\partial t} x = k + \left( \frac{p-1}{r} + \frac{r}{2} \right) \sin \theta - \left( \frac{u}{2} + \frac{q-1}{u} \right) \cos \theta, \qquad x \in \gamma_t,$$

which is nonsingular away from the axes (see (8)).

Let $\sigma$ be the ray $u = r \tan \alpha_s$, so $C_{p,q} = M(\sigma)$.

**Proof.** Let $\{M_t\}_{t \in \mathbf{R}}$ be one of the shrinking-expanding solutions constructed above, with $M_t = M(\gamma_t)$. We cannot apply Angenent's Theorem directly, since $\gamma_0$ meets $\sigma$ at the origin, where both the curves and the equation are singular.

Recall the wiggly minimal surface $Q = M(v)$ of section B, and write $\varepsilon = \text{graph}(v)$. Then $\#(\gamma_t \cap \varepsilon)$ is nonincreasing for $t \neq 0$. By passing to limits near $t = 0$ we get

$$\#(\gamma_s \cap \varepsilon) \leq \#(\gamma_0 \cap \varepsilon) \leq \#(\gamma_t \cap \varepsilon), \qquad t < 0 < s.$$

In fact, we have

$$\lim_{s \to \infty} \#(\gamma_s \cap \varepsilon) + 2 \leq \lim_{t \to -\infty} \#(\gamma_t \cap \varepsilon)$$

since $\gamma_t$ is tangent to $\varepsilon$ once coming down and once coming back up, at the point $r = 0$, $u = v(0)$. (The strict decrease can be verified by hand, without having to extend Angenent's theorem to the singular point $r = 0$.)

This proves the proposition, because

$$\lim_{s \to \infty} \#(\gamma_s \cap \varepsilon) = \#(\gamma_1 \cap \sigma), \qquad \lim_{t \to -\infty} \#(\gamma_t \cap \varepsilon) = \#(\gamma_{-1} \cap \sigma).$$

In fact, the number of intersections drops by at least two.

**Remark.** 1. As remarked above, there are also expanding evolutions of the form $r = r(u)$; in this case, there is only one moment of tangency, so the number of intersections drops by at least one.

2. We would like a more exact estimate of the number of continuations. In light of the above results, one might conjecture (a) for each $k \geq 2$, there is exactly one self-shrinker $N_t^k$ of the form $M(\{u = u(r)\})$ which is asymptotic to a cone at infinity and has $k$ intersections with the cone, and (b) $N_t^k$ has exactly $2k-1$ continuations $P_t$, with the following numbers of intersections with the cone:

$$0, 2, 2, \ldots, k-2, k-2, \qquad \text{pulling one way,}$$

$$1, 1, 3, 3, \ldots, k-1, k-1, \qquad \text{pulling the other way.}$$

Finally we show that the cone $C_{p,q}$ itself can never arise as a blowup.

**Proposition.** *Suppose $\{M_t\}_{t\in[-b,0)}$ is a smooth, compact, rotationally symmetric mean curvature flow that becomes singular at $t=0$. Then $M_t^\lambda$ cannot converge to the cone $C_{p,q}$.*

**Proof.** The basic idea, due to Grayson, is this: the self-shrinking solutions we have constructed "protect" the cone by having too many wiggles.

Let $M_t = M(\gamma_t)$. For $t < 0$, $\gamma_t$ meets the axes orthogonally (or not at all) and never contains $(0,0)$.

By Angenent [A1], $\gamma_t \cap \sigma$ is finite and and the intersection is transverse for almost every $t \in (-b, 0)$. Suppose in particular that this holds at $t = t_0$.

For large $k$, $\delta_k$ is close to $\sigma$ in $C^1_{loc}$ and $u/r$ is close to $\lambda_s$ in $C^1(\mathbf{R}^2 \setminus B_1)$. It follows that for $k$ sufficiently large,

$$\#(\gamma_{t_0} \cap \delta_k) = \#(\gamma_{t_0} \cap \sigma) \ll \#(\delta_k \cap \sigma).$$

Since $\gamma_t$ and $\sqrt{-t} \cdot \delta_k$ both solve (19), by Angenent [A1], $\#(\gamma_t \cap \sqrt{-t} \cdot \delta_k)$ is nonincreasing. On the other hand, if $(-t)^{-1/2} \cdot \gamma_t$ converges to $\sigma$, then $\#((-t)^{-1/2} \cdot \gamma_t \cap \delta_k) \to \#(\sigma \cap \delta_k)$, which is a contradiction. This proves the proposition.

A similar proof, using the wiggly curve $\varepsilon$ instead of $\delta_k$, shows that if $M_t$ (possibly noncompact) is smooth and rotationally symmetric for $t < 0$, then

$$M_t \not\to C_{p,q} \qquad \text{as } t \to 0^-$$

provided that $M_t$ originally has a finite number of intersections with the cone.

**Example.** Interestingly, we can attain $C_{p,q}$ exactly if we allow triple junctions, in the case $p = q$ and $n = 4$ or $6$, with thanks to J. Sullivan and R. Kusner. Consider the self shrinkers

$$u_{\sqrt{n-2}}(r) = \sqrt{n-2} \qquad \text{(Cylinder)}$$
$$u_{\sqrt{2(n-1)}}(r) = \sqrt{2(n-1) - r^2} \qquad \text{(Sphere)}$$

which cross $u = r$ at $135°$ and $90°$ respectively. Somewhere between, there is a curve $u_a$ that hits $u = r$ at $120°$, at some point $(b, b)$. Let $\gamma_1$ denote the arc of $u_a$ from $(0, a)$ to $(b, b)$, let $\gamma_2$ be its reflection across $u = r$, and let $\gamma_3$ be the part of $u = r$ from $(b, b)$ to infinity. Then

$$N_t = \sqrt{t} \cdot M(\gamma_1 \cap \gamma_2 \cap \gamma_3), \qquad t < 0,$$

is a self-shrinker with triple junctions that meet at $120°$, such that $N_t \to C_{p,q}$ as $t \to 0$. A trapped bubble shrinks to nothing at the origin, and the cone explodes in infinitely many possibilities.

Tom Ilmanen

Evanston

August, 1995.